\documentclass[a4paper]{article}

\usepackage[english]{babel}
\usepackage[utf8x]{inputenc}
\usepackage[T1]{fontenc}

\usepackage{graphics}

\usepackage{algorithm}
\usepackage{algpseudocode}

\algnewcommand{\Multcomment}[1]{\\\hfill \begin{minipage}[t]{0.8\textwidth}/* {\itshape #1} */\end{minipage}}

\algnewcommand{\algorithmicgoto}{\textbf{go to}}%
\algnewcommand{\Goto}[1]{\algorithmicgoto~\ref{#1}}%

\usepackage{enumitem} 
\usepackage{amsmath}
\usepackage{amssymb}
\usepackage{graphicx}
\usepackage[colorinlistoftodos]{todonotes}
\usepackage[colorlinks=true, allcolors=blue]{hyperref}
\usepackage{color}
\usepackage{multirow}
\usepackage[round, sort,comma,authoryear]{natbib}
\usepackage{siunitx} 
\usepackage{verbatim}
\usepackage{tikz}

\newcommand{\footnoteremember}[2]{
\footnote{#2}
  \newcounter{#1}
  \setcounter{#1}{\value{footnote}}
}
\newcommand{\footnoterecall}[1]{
\footnotemark[\value{#1}]
}

\title{Fast Continuous and Integer L-shaped Heuristics Through Supervised Learning}

\author{
Eric Larsen\footnoteremember{cirrelt}{Department of Computer Science and Operations Research and CIRRELT, Universit\'e de Montr\'eal}
\and
Emma Frejinger\footnoterecall{cirrelt}\footnote{Corresponding author. Email: frejinger.umontreal@gmail.com}
\and
Bernard Gendron\footnoterecall{cirrelt}
\and
Andrea Lodi\footnoteremember{Lodi}{CERC, Polytechnique Montr\'eal and Jacobs Technion-Cornell Institute, Cornell Tech and Technion}}

\date{\today}

\begin{document}

\maketitle

\begin{abstract}
We propose a methodology at the nexus of operations research and machine learning (ML) leveraging generic approximators available from ML to accelerate the solution of mixed-integer linear two-stage stochastic programs. We aim at solving problems where the second stage is highly demanding. Our core idea is to gain large reductions in online solution time while incurring small reductions in first-stage solution accuracy by substituting the exact second-stage solutions with fast, yet accurate supervised ML predictions. This upfront investment in ML would be justified when similar problems are solved repeatedly over time, for example, in transport planning related to fleet management, routing and container yard management.

Our numerical results focus on the problem class seminally addressed with the integer and continuous L-shaped cuts. Our extensive empirical analysis is grounded in standardized families of problems derived from \emph{stochastic server location} (SSLP) and \emph{stochastic multi knapsack} (SMKP) problems available in the literature. The proposed method can solve the hardest instances of SSLP in less than 9\% of the time it takes the state-of-the-art exact method, and in the case of SMKP the same figure is 20\%. Average optimality gaps are in most cases less than 0.1\%.
\end{abstract}


\section{Introduction} \label{sec:intro}
Decision-making problems in a broad range of application domains are subject to uncertainty. Prominent examples are supply chain management and transport planning that during the COVID-19 pandemic received media attention for their widespread disruptions and lack of resiliency. Real-life instances of related decision-making problems, such as fleet management and network design problems, are known to be hard to solve, even in a deterministic form. Stochastic formulations are therefore rarely used in practice, despite their well-documented advantages when there is uncertainty in problem parameters. 

We consider a general stochastic linear two-stage problem (P) of the form \citep[notation borrowed from][]{AnguEtAl2016}
\begin{alignat}{2} \nonumber
\min_{x, z, \theta}~ &  \lbrace cx + dz + \theta\rbrace \\ 
\label{eq:stage1_const1}   \text{s.t. } & Ax + Cz \leq b, \\
\label{eq:stage1_const2}  & Q(x) - \theta \leq 0, \\
\label{eq:stage1_const3} & x \in \lbrace 0, 1 \rbrace^n, \\
\label{eq:stage1_const4} & z \geq 0, ~~ z \in \mathcal{Z},
\end{alignat} 
where the second-stage subproblem (S) is
\begin{equation}
Q(x) :\equiv \mathbb{E}_{\xi} [ \min_y \lbrace q_{\xi} y :
 W_{\xi}y \geq h_{\xi} - T_{\xi} x, y \in \mathcal{Y} \rbrace ],
\label{eq:Benders_subproblem}
\end{equation}
and $\mathcal{Z}$ and $\mathcal{Y}$ may embody integrality constraints and a number of additional constraints pertaining to $z$ and $y$, respectively. The coupling constraints $W_{\xi}y \geq h_{\xi} - T_{\xi}x$ tie the first and second stages. We concentrate on formulations where (i) stochastic coefficients $\xi \equiv: \{q_{\xi}, W_{\xi},T_{\xi},h_{\xi}\}$ with finite support can occur in the second-stage objective or constraints, (ii) both stages can feature integer or continuous variables, (iii) the variables coupling the first and second stage are binary-valued, and (iv) both (S) and its primal continuous relaxation (RS)
\begin{equation}
\widetilde{Q}(x) :\equiv \mathbb{E}_{\xi} [ \min_y \lbrace q_{\xi} y : W_{\xi}y \geq h_{\xi} - T_{\xi} x, y \in \widetilde{\mathcal{Y}} \rbrace ]
\label{eq:Benders_relaxed_subproblem}
\end{equation}
feature relatively complete recourse ($\widetilde{\mathcal{Y}}$ is equal to $\mathcal{Y}$ except for the removal of every integrality constraint).

Common to many real-world settings is the solution of instances of (P) sharing similar patterns. Consider for example the tactical-operational planning of transport systems. This similarity lends itself to machine learning (ML) where data on the distribution of instances of interest could be leveraged offline to gain online speed-up. In this context, our goal is to devise a general ML-based matheuristic to solve mixed-integer linear two-stage stochastic programs, where the second stage is highly demanding computationally. More precisely, our core idea is to gain large reductions in online solution time while incurring small reductions in first-stage solution accuracy by substituting the solutions of the latter stages of Benders decomposition with fast, yet accurate approximations arising from supervised ML. 

\subsection{Relations with existing literature}

Let us turn our attention (i) to exact methods currently available in the operations research (OR) literature addressing problem class (P) and to heuristic methods from the (ii) OR and the (iii) ML literatures and to the relations of our own heuristic method therewith.

\subsubsection{Exact methods}
Problem (P) is seminally addressed in \cite{LapoLouv1993}  with the integer L-shaped method operating on the Benders decomposition of (P). Let $L$ denote a lower bound on $Q(x)$. The master (first-stage) problem (M) is then
\begin{alignat}{2} \nonumber
\min_{x, z, \theta}~ & \lbrace cx + dz + \theta\rbrace \\
\nonumber \text{s.t. } & (\ref{eq:stage1_const1}), (\ref{eq:stage1_const3}), (\ref{eq:stage1_const4}), \\
&\label{eq:optimality_cuts} \Pi x - \mathbf{1} \theta \leq \pi_0, \\
&\theta \geq L,
\end{alignat}
where $\mathbf{1}$ denotes a column vector of ones. The set of constraints (\ref{eq:optimality_cuts}) is initially (i.e., at the root node of the algorithm) empty (unless initial cuts have been added, for instance, in a multi-phase solution process), and progressively populated with optimality cuts as a branch-and-Benders-cut process advances. \cite{LapoLouv1993} apply the integer L-shaped optimality cut at a candidate solution $x^*\in \{0,1\}^n$
\begin{equation}
\left(Q(x^*) - L \right) \Big(\sum_{i \in S(x^*)} x_i - \sum_{i \notin S(x^*)} x_i - |S(x^*)| \Big) + Q(x^*) \leq \theta,
\label{eq:integer_L_shaped_cut}
\end{equation}
where $S(x^*) :\equiv \lbrace i: x_i^* = 1 \rbrace$ and $Q(x^*)$ is the optimal value of (S) at $x^*$, whenever $x^*$ turns out to be an invalid first-stage solution.

To the best of our knowledge, \cite{AnguEtAl2016} is currently the most up-to-date and complete overview available about the application of Benders decomposition methods to problem class (P). They propose an algorithm that alternates between continuous and integer L-shaped cuts. At a given candidate integral first-stage solution $x^*$ that turns out to be invalid, their alternating cuts strategy first computes a subgradient cut derived from (RS). In particular, if $\mathcal{Y}$ only imposes that second-stage variables in subproblem (S) are binary, then a continuous L-shaped optimality mono-cut \citep[pp.~183-184]{BirgLouv2011} is defined by
\begin{equation}
\mathbb{E}_{\xi} [\phi_{\xi}(h_{\xi} - T_{\xi} x) - \mathbf{1}' \psi_{\xi}] \leq \theta,
\label{eq:continuous_L_shaped_cut}
\end{equation}
where $\phi_{\xi}, \psi_{\xi}$ are the solutions to the duals (DRS) corresponding to RS, evaluated at $x^*$,
\begin{equation}
\max_{\phi_{\xi}, \psi_{\xi}} \lbrace \phi_{\xi} (h_{\xi} - T_{\xi} x^*) - \mathbf{1}' \psi_{\xi} : \phi_{\xi} W_{\xi} - \psi_{\xi} \leq q_{\xi}, \phi_{\xi} \geq 0, \psi_{\xi} \geq 0 \rbrace,
\label{eq:Benders_relaxed_subproblem_dual}
\end{equation}
$\forall \xi$. Here, the dual variables $\psi_{\xi}, \forall \xi$ are attached to the upper bound (equal to 1) embodied in $\widetilde{\mathcal{Y}}$ and constraining the second-stage variables in (RS).

The alternating algorithm of \cite{AnguEtAl2016} resorts to the calculation of the integer L-shaped cut (\ref{eq:integer_L_shaped_cut}) (which is then imposed in the first-stage problem) only if the continuous L-shaped mono-cut (\ref{eq:continuous_L_shaped_cut}) fails to separate the invalid first-stage solution. Otherwise, the continuous L-shaped cut is added to the first stage. This algorithm is generally faster than the standard integer L-shaped algorithm since (i) calculating the solution of the integral second-stage problem usually requires considerably more time than calculating the solution of its continuous relaxation and (ii) a binding continuous L-shaped cut will generally create a deeper separation than the corresponding integer L-shaped cut. 

\cite{AnguEtAl2016} also propose a new type of integer optimality cut, based on disjunctive programming, that can be used in place of (\ref{eq:integer_L_shaped_cut}), either in the integer L-shaped method or in their alternating algorithm. However, their experimental results indicate that the introduction of these cuts has a marginal or relatively small impact on computing times, depending on the class of problems considered.

\subsubsection{Heuristic methods from OR}

We turn our attention to heuristic algorithms that have been proposed to solve (P). A strong alternative to the aforementioned exact methods from the OR literature is the \textit{progressive hedging} (PH) algorithm \citep{RockWets91, WatsWood11} in view of its broad applicability, the interest it has attracted, its advanced and continuing development, and of the availability of an extensive computational platform. The \emph{dual decomposition} algorithm \citep{CaroRudi99} is also a possible contender. However, the documented performance of PH, provided it is favorably tuned, apparently surpasses that of dual decomposition for certain classes of problems, notably that of \emph{stochastic server location problems} on which we draw as a test bench in our experiments \citep{TorrEtAl22, Wood16}.

\paragraph{Commonalities of the proposed method in OR.}
In OR, the idea of approximating all or part of a complex programming problem or of its solution has chiefly been investigated in simulation optimization and, to a much lesser extent, in bilevel programming. In our method, an ML predictor plays a role comparable to the objects alternatively known as \emph{metamodels}, \emph{surrogate models} or \emph{response surfaces} \citep{BartMeck2006, Bart2020} that are encountered in simulation optimization \citep{AmarEtAl2016}. Likewise, our ML predictor describes a map between the inputs and output(s) of an optimization process whose exact form is fitted to data. In contrast, we are neither performing optimization over the ML predictor itself as a surrogate for the map between inputs and output(s) of the original optimization problem nor inserting it as a surrogate constraint as a stand in for a subproblem occurring in the original optimization problem. It is used, instead, as an auxiliary instrument to accelerate computations by generating at a high speed the intermediate inputs that are required in the broader optimization process governed by a Benders decomposition. This use as an auxiliary instrument is exploited in \citet{LarsEtAl2022a} and \citet{LarsFrej2019}. From this perspective, and to our knowledge, our closest neighbors are \citet{SinhEtAl2016} and \citet{SinhEtAl2020}, where metamodels approximate the lower-level reaction set mapping and the lower-level value function in evolutionary algorithms handling deterministic, bilevel programming problems featuring continuous variables.

\subsubsection{Heuristic methods from ML}

There is an abundant and steadily growing literature in ML about constrained optimization, and combinatorial optimization in particular. It is surveyed in \citet{BengEtAl2021,  VesseEtAl2020, KotaEtAl2021a,  KotaEtAl2021b}. \citet{BengEtAl2021} and \citet{KotaEtAl2021a} distinguish two main areas in which ML operates. The first one, \emph{ML alongside optimization algorithms} \citep{BengEtAl2021}, i.e., \emph{ML-augmented constrained optimization} \citep{KotaEtAl2021a}, designates the utilization of ML to inform the decisions made by optimization algorithms in order to enhance their efficiency. Our contribution is not quite of this nature. The second main area is \emph{end-to-end learning}, whose defining characteristic is simply the application of ML within calculations towards the solution of optimization problems. This better corresponds to our matheuristic algorithm: we utilize ML to approximate a core component of the Benders decomposition, that is the solution of second stage, given first-stage coupling variables, to achieve an approximate first-stage solution for the stochastic programming problem at hand.

\citet{KotaEtAl2021a} further distinguish two subareas in end-to-end constrained optimization learning: (i) \emph{predicting solutions} and (ii) \emph{predict-and-optimize}. On the one side, (i) designates the application of supervised or reinforcement learning (RL) to predict the overall solutions of constrained optimization problems. \citet{CappEtAl2020, MazyEtAl2021, WangTang2021} detail the application of RL in this context. Graphical neural networks (GNN) are especially well-adapted to the representation of combinatorial problems and their application is detailed in \citet{CappEtAl2021, PengEtAl2021, SchuEtAl2022}. On the other side, (ii) designates the search for solutions to partially defined problems while accounting jointly for their unknown aspects or parameters. Predecessors of methods currently in (ii) belong to an area known as \emph{predict-then-optimize}, where prediction of the unknown aspects of the optimization problem and optimization itself are performed sequentially. Recent information on (ii) is available in \citet{SPTO, Fior2022, ZhanEtAl2022}. Whereas methods in subarea (ii) involve the application of mathematical programming solvers, those in (i) rely entirely on optimization of ML architectures. It is worth noting that an assignment of papers to specific subareas is often not totally precise, because some works can be interpreted differently and might fit into multiple categories.

The following recent publications lie at the intersection of stochastic programming and ML. Given a set of scenarios, \citet{BengEtAl2020} propose to predict a representative scenario to form a smaller surrogate problem that is easier to solve with a general purpose MIP solver. The definition of a representative scenario depends on the structure of the problem. This approach can be assigned to both ML-augmented constrained optimization and end-to-end learning. In constrast, \citet{NairEtAl2018} and \citet{DaiEtAl2021} belong  to subarea (i) of end-to-end learning. \citet{NairEtAl2018} propose an RL algorithm for two-stage stochastic programs with unconstrained binary decisions and hence does not address the general problem class (P). \citet{DaiEtAl2021} address multi-stage stochastic programs. Although this class is a generalization of (P), the respective contexts and solution approaches are distant: The authors exploit the commonalities between the problems successively solved by the stochastic dual dynamic programming algorithm to progressively train an efficient approximator of the value function. In a very recent paper, \citet{DumoEtAl2022} propose to approximate the second-stage solution value with a feed-forward neural network \citep[in a fashion similar to][]{LarsEtAl2022a}. Using the results in \cite{FiscJo2018}, they represent the feed-forward network as a mixed-integer programming (MIP) problem, which allows them to integrate it in the MIP formulation of the first-stage problem.

\paragraph{Relation of proposed method with ML.}

The method we propose falls outside of subareas (i) and (ii) currently explored in ML. Since it involves an ML predictor whose role is to inform an overall solution process handled by a solver, it is somewhat more closely connected to (ii). Yet, the uncertainty of the problems it addresses is fully specified rather than discovered through prediction and, instead of missing data, our method predicts intermediate values (e.g., expectations) that are instrumental to the overall optimization process governed by the underlying mathematical programming algorithm. It also optimizes the predictor and the overall solution sequentially, rather than jointly. To the best of our knowledge, \citet{DumoEtAl2022} is the only other work addressing the general problem class (P) with ML. Beyond approximating the second-stage solution value with ML, the differences between their approach and the methodology we propose run deep.

\subsection{Contributions}
We propose an ML-based matheuristic version of the L-shaped method that we call \emph{ML-L-Shaped}. It is designed to effectively solve problems of class (P) featuring computationally demanding second-stage problems. More precisely, ML-L-Shaped is a heuristic version of both the standard integer L-shaped method \citep{LapoLouv1993} and the algorithm with alternating continuous and integer L-shaped cuts \citep{AnguEtAl2016}, where we replace the costly computations, that is, solution values of (S) -- $Q(x)$ -- and (RS) -- $\widetilde Q(x)$ -- by fast ML predictions, and solution values of (DRS) -- $\psi_{\xi}, \phi_{\xi}, \forall \xi$ -- by fast ML predictions of low-dimensional reductions.

Focusing on input-output structure and dimension reduction, we describe how to generate labeled data about the second-stage problems so as to use standard supervised ML. We also discuss how supervised ML can utilize a method presented in  \cite{LarsEtAl2022a} to reduce the cost of generating labeled data.

Our numerical analyses are grounded in the two standard classes of problems that are examined in \cite{AnguEtAl2016}: the \textit{stochastic server location problems} (SSLP) and the \textit{stochastic multi knapsack problems} (SMKP). We compare the performance of ML-L-Shaped to the exact counterpart as well as to a PH algorithm. Our method can solve the hardest instances of SSLP in less than 9\% of the time it takes the state-of-the-art exact method, and in the case of SMKP the same figure is 20\%. Average optimality gaps are in most cases less than 0.1\%.

\subsection{Structure of the paper}

The remainder of the paper is structured as follows: Section~\ref{sec:decomposition_algorithms} introduces the ML-L-Shaped method. Section~\ref{sec:experimental_setting} describes how we operationalize the method, i.e., the experimental setting: We review the standardized problem families on which the detailed analysis is based, the generation of data from each family and the input-output structures of the ML approximators.
Section~\ref{sec:ml_approximation} outlines the training and validation of ML predictors for the solutions of integral and continuously-relaxed second-stage problems. Section~\ref{sec:detailed_analysis} details and compares the experimental results achieved by the exact and ML-assisted versions of the algorithms presented in Section~\ref{sec:decomposition_algorithms}. We compare their performances to that of an alternative heuristic algorithm (PH) in Section~\ref{sec:comparison_with_alternative_heuristic}. Section~\ref{sec:Conclusion} summarizes our findings and outlines possible extensions and implementations of our core idea.

\section{ML-based L-shaped method}
\label{sec:decomposition_algorithms}

This section introduces our ML-L-Shaped method. We distinguish between the heuristic versions of the standard (Std) and alternating (Alt) cut strategies in the argument of the main procedure of the branch-and-Benders-cut process (Algorithm~\ref{alg:benders}) by setting $\textit{isAlt}$ to \emph{true} in the case of alternating cuts, and to \emph{false} otherwise. Those strategies are standard and essentially correspond to the ones of \cite{AnguEtAl2016}, except for the heuristic callback (Algorithm~\ref{alg:heuristic_callback_v2}) introduced in Step~\ref{a:callback}. 

\begin{algorithm}
\caption{Benders decomposition: Main} \label{alg:benders}
\begin{algorithmic}[1]
\Procedure{Main}{$\boldsymbol{isAlt}, \boldsymbol{\mu}, \boldsymbol{\nu}$}

\State Compute or retrieve the lower bound $L$ for the objective value of (P).
\State \parbox[t]{313pt}{Initialize a branch-and-cut process with a global node tree for (M). This creates the repository of leaf nodes, say $R$. The latter initially contains only the root node.\strut} \label{a:main_initialize} 

\State $\textit{UB} \gets \infty$ \Comment{First-stage upper bound}

\State $(x^{**}, z^{**}) \gets \varnothing$ \Comment{First-stage incumbent solution}

\If{$R = \emptyset$} \label{select}
\State \Goto{final}
\Else
\State{Select a node from $R$.}
\EndIf

\State  \parbox[t]{310pt}{Compute the current optimal solution $(x^*, \theta^*)$ to (M) for the node at hand. \strut} \label{SolM} 

\If{ $(cx^* + dz^* + \theta^*) \geq \textit{UB}$}
\State Discard the node from $R$.
\State \Goto{select} 
\EndIf
\If{$(x^*, z^*)$ is not integral} \label{alg:not_integral}
\State  \parbox[t]{290pt}{Partition the domain of $(x, z)$ in (M) or add MIP-based cuts. Accordingly add newly defined nodes to $R$ or update existing nodes in $R$. \strut}
\State \Goto{select}
\EndIf

\State \Call{HeuristicCallback}{$\boldsymbol{isAlt}, \boldsymbol{\mu}, \boldsymbol{\nu}$} \label{a:callback} 

\State \Goto{select}
\State  \parbox[t]{300pt}{Retrieve the final first-stage incumbent solution $(x^{**}, z^{**})$. Compute the final overall value $cx^{**} + dz^{**} + Q(x^{**})$. \strut} \label{final}

\EndProcedure
\end{algorithmic}
\end{algorithm}

Let $Q^{\textit{ML}}(x^*)$ denote the value of $Q(x^*)$ predicted by ML. We define a heuristic version of the exact integer L-shaped cut (\ref{eq:integer_L_shaped_cut}) as
\begin{equation}
\left( Q^{\textit{ML}}(x^*) - L \right) \Big(\sum_{i \in S(x^*)} x_i - \sum_{i \notin S(x^*)} x_i - |S(x^*)| \Big) + Q^{\textit{ML}}(x^*) \leq \theta.
\label{eq:integer_L_shaped_cut_heur}
\end{equation}

Similarly, let $\widetilde{Q}^{\textit{ML}}(x^*)$ denote the value of $\widetilde{Q}(x^*)$ predicted by ML
in (\ref{eq:Benders_relaxed_subproblem}). Now, we might consider defining a heuristic version of the exact continuous L-shaped mono-cut (\ref{eq:continuous_L_shaped_cut}) as
\begin{equation}
\mathbb{E}_{\xi} [\phi^{\textit{ML}}_{\xi}(h_{\xi} - T_{\xi} x) - \mathbf{1}' \psi^{\textit{ML}}_{\xi}] \leq \theta,
\end{equation}
 where $\phi^{\textit{ML}}_{\xi}$ and $\psi^{\textit{ML}}_{\xi}$ would denote predictions of $\phi_{\xi}$ and $\psi_{\xi}$, respectively, appearing in (\ref{eq:Benders_relaxed_subproblem_dual}), $\forall \xi$.
 
 However, it is advantageous to note that detailed knowledge of the potentially high-dimensional $\phi_{\xi}$ and $\psi_{\xi}$, $\forall \xi$ is not necessary to construct the mono-cut  (\ref{eq:continuous_L_shaped_cut}). The three reductions $\mathbb{E}_{\xi} [\phi_{\xi}h_{\xi}]$, $\mathbb{E}_{\xi} [\phi_{\xi}T_{\xi}]$ and $\mathbb{E}_{\xi}  [\mathbf{1}' \psi_{\xi}]$ suffice. Furthermore, whenever $h_{\xi}$  or $T_{\xi}$ are non-stochastic, that is $h_{\xi} \equiv h, \forall \xi$ or $T_{\xi} \equiv T, \forall \xi$, it is sufficient to calculate $\mathbb{E}_{\xi} [\phi_{\xi}]$ and $\mathbb{E}_{\xi} [\psi_{\xi}]$, respectively. By exploiting these reductions, the task of building an ML approximation for (\ref{eq:continuous_L_shaped_cut}) can be eased. Hence, we define the following heuristic version of the exact continuous L-shaped mono-cut (\ref{eq:continuous_L_shaped_cut}) as
\begin{equation}
\lbrace \mathbb{E}_{\xi} [\phi_{\xi}h_{\xi}]\rbrace ^{ML} - \lbrace \mathbb{E}_{\xi} [\phi_{\xi}T_{\xi}]\rbrace ^{ML} x - \lbrace \mathbb{E}_{\xi} [\mathbf{1}' \psi_{\xi}]\rbrace ^{ML} \leq \theta,
\label{eq:continuous_L_shaped_cut_heur}
\end{equation}
where $\lbrace \mathbb{E}_{\xi} [\phi_{\xi}h_{\xi}]\rbrace ^{ML}$, $\lbrace \mathbb{E}_{\xi} [\phi_{\xi}T_{\xi}]\rbrace ^{ML}$ and $\lbrace \mathbb{E}_{\xi} [\mathbf{1}' \psi_{\xi}]\rbrace ^{ML}$ are ML predictions of the reductions $\mathbb{E}_{\xi} [\phi_{\xi}h_{\xi}]$, $\mathbb{E}_{\xi} [\phi_{\xi}T_{\xi}]$ and $\mathbb{E}_{\xi} [\mathbf{1}' \psi_{\xi}]$, respectively.

Alternatively, whenever $h_{\xi}$  and $T_{\xi}$ are non-stochastic, that is $h_{\xi} \equiv h, \forall \xi$ or $T_{\xi} \equiv T, \forall \xi$, simplifications ensue and we calculate the heuristic version of the exact continuous L-shaped mono-cut (\ref{eq:continuous_L_shaped_cut}) as 
\begin{equation}
\lbrace \mathbb{E}_{\xi} [\phi_{\xi}]\rbrace ^{ML} (h - T x) - \lbrace \mathbb{E}_{\xi} [\mathbf{1}' \psi_{\xi}]\rbrace ^{ML} \leq \theta,
\label{eq:continuous_L_shaped_cut_heur_nonstoch_h_T}
\end{equation}
where $\lbrace \mathbb{E}_{\xi} [\phi_{\xi}]\rbrace ^{ML}$ and $\lbrace \mathbb{E}_{\xi} [\mathbf{1}' \psi_{\xi}]\rbrace ^{ML}$ are ML predictions of the reductions $\mathbb{E}_{\xi} [\phi_{\xi}]$ and $\mathbb{E}_{\xi} [\mathbf{1}' \psi_{\xi}]$, respectively. These simplifications occur in the application reported in Section~\ref{families_of_smkp_instances}.

Since the predictions occurring at Step~\ref{a:contPred} of Algorithm~\ref{alg:heuristic_callback_v2} may be yielded by a single prediction model, we compute all of them simultaneously. In practice, at Steps~\ref{a:contPred}~and~\ref{heur_is_not_oppor2} of Algorithm~\ref{alg:heuristic_callback_v2}, C-language bindings launch GPU computations returning ML predictions.

In addition to the argument determining the variant of the algorithm (standard integer L-shaped or alternating cuts strategy), there are two hyperparameters: $0 \ll \mu \leq 1$ and $0 \ll \nu \leq 1$. They control the likelihood of introducing incorrect integer L-shaped or continuous L-shaped cuts. This could occur even if the predictions yielded by ML are of high quality and feature low average absolute relative errors. Concretely, since we are solving minimization problems, we wish to control the likelihood that $\widetilde Q^{\textit{ML}}(x^*)$ and $Q^{\textit{ML}}(x^*)$ overestimate the corresponding exact values. We achieve this by uniformly shifting down the ML predictions. The exact values of these shifts specified with $\mu$ and $\nu$ are estimated from a preliminary tuning process on an independent set of problem instances, starting with a value of one. Note that we obtain the exact counterparts of the heuristic algorithm by fixing $\mu=\nu=1$ and by using exact values $\widetilde{Q}(x^*)$, $\mathbb{E}_{\xi} [\phi_{\xi}h_{\xi}]$, $\mathbb{E}_{\xi} [\phi_{\xi}T_{\xi}]$, $\mathbb{E}_{\xi} [\mathbf{1}' \psi_{\xi}]$, $\mathbb{E}_{\xi} [\phi_{\xi}]$, $\mathbb{E}_{\xi} [\mathbf{1}' \psi_{\xi}]$, $Q(x^*)$ instead of predictions, as well as cuts (\ref{eq:continuous_L_shaped_cut}) and (\ref{eq:integer_L_shaped_cut}) in Steps~\ref{alg:heurContL} and~\ref{alg:heurIntL} of Algorithm~\ref{alg:heuristic_callback_v2}, respectively.

\begin{algorithm}
\caption{Benders decomposition: Heuristic callback} \label{alg:heuristic_callback_v2}
\begin{algorithmic}[1]
\Procedure{HeuristicCallback}{$\boldsymbol{isAlt}, \boldsymbol{\mu}, \boldsymbol{\nu}$}

\If{$!\boldsymbol{isAlt}$} 
\State \Goto{heur_is_not_oppor2} 
\EndIf
\State \parbox[t]{290pt}{\text{Compute predictions} \Comment{Alternating cut strategy} \\
$\widetilde{Q}^{\textit{ML}}(x^*)$, $\lbrace \mathbb{E}_{\xi} [\phi_{\xi}h_{\xi}]\rbrace ^{ML}$, $\lbrace \mathbb{E}_{\xi} [\phi_{\xi}T_{\xi}]\rbrace ^{ML}$, $\lbrace \mathbb{E}_{\xi} [\mathbf{1}' \psi_{\xi}]\rbrace ^{ML}$ \\
or \\
$\widetilde{Q}^{\textit{ML}}(x^*)$, $\lbrace \mathbb{E}_{\xi} [\phi_{\xi}]\rbrace ^{ML}$, $\lbrace \mathbb{E}_{\xi} [\mathbf{1}' \psi_{\xi}]\rbrace ^{ML}$.}
\label{a:contPred}   
\If{$\nu \widetilde{Q}^{\textit{ML}}(x^*) > \theta^*$} 
\State \parbox[t]{290pt}{Add a heuristic continuous L-shaped mono-cut (\ref{eq:continuous_L_shaped_cut_heur}) or (\ref{eq:continuous_L_shaped_cut_heur_nonstoch_h_T}). \strut} \label{alg:heurContL}
\State \Return
\EndIf
\State \text{Compute prediction $Q^{\textit{ML}}(x^*)$}. \label{heur_is_not_oppor2} \Comment{Integer L-shaped method}
\If {$\mu Q^{ML}(x^*) \leq \theta^*$} 
\If {$cx^{*} + dz^{*} + \theta^* < \text{\it UB}$}
\label{alg:checkUB}
\State $\text{\it UB} \gets cx^{*} + dz^{*} + \theta^*$ \Comment{Update upper bound} \label{alg:UpdateUB}
\State $(x^{**}, z^{**}) \gets (x^*, z^*)$ \Comment{Update incumbent solution}
\EndIf
\Else
\State \text{Add a heuristic integer L-shaped cut (\ref{eq:integer_L_shaped_cut_heur}).} \label{alg:heurIntL}
\EndIf
\EndProcedure
\end{algorithmic}
\end{algorithm}

Algorithm~\ref{alg:benders} terminates with a feasible solution, provided that Step~\ref{alg:UpdateUB} of Algorithm~\ref{alg:heuristic_callback_v2} is reached at least once. The likelihood of reaching that step is controlled by the aforementioned shift coefficients. In the unlikely event that the algorithm terminates without a feasible solution, it can be reapplied, after decreasing the values of $\mu$ and $\nu$. This occurred twice over all problem instances considered in Section~\ref{sec:detailed_analysis}.

We now turn to the differences between the heuristic integer L-shaped and alternating cut algorithms. First, recall that the alternating cut strategy was proposed by \cite{AnguEtAl2016} to avoid costly computations of $Q(x^*)$. On the contrary, predictions $Q^{\textit{ML}}(x^*)$ are fast to compute (in the order of a few milliseconds). Moreover, as we further discuss in Section~\ref{sec:experimental_setting}, the task of predicting $Q(x^*)$ is easier than that of predicting $\widetilde{Q}(x^*)$ and the reductions of $\phi_{\xi}, \psi_{\xi}, \forall \xi$. This a priori favors the heuristic version of the standard integer L-shaped method over the alternating cuts strategy. However, the latter is likely to be useful when the first-stage problem is hard. Indeed, integer L-shaped cuts suppress only one first-stage solution at a time whereas continuous L-shaped cuts can be stronger. This can make a difference when the number of heuristic integer L-shaped cuts is very large if used alone. In these circumstances, the solution process can be hindered by their large number despite the very high speed of the individual computations.

We consider two further variants of the proposed algorithm. Both involve two phases, where Algorithm~\ref{alg:benders} is used in a first phase to produce a feasible solution. In the first variant, this solution is used alone to warm start the exact standard integer L-shaped method, or the one with alternating cuts. The resulting approach is then exact. In the second variant, in addition to supplying a warm-start solution, we introduce a probabilistic lower bound on the value of the first-stage objective in the exact solution process. We obtain this probabilistic lower bound for a given problem family by computing the empirical distribution of exact objective values from a preliminary, independent set of instances and calculating a 10\% one-sided Chebyshev lower confidence bound \citep{Ngo2011}. Such confidence bounds are pessimistic and the resulting solutions can be interpreted as having a probability of at least 90\% of being exact.

\section{Experimental setting: problem families, data generation and input-output structure of ML predictors}
\label{sec:experimental_setting}

We ground our numerical analysis of the ML-L-Shaped method in the standard problem sets of classes SSLP and SMKP that are examined in \cite{AnguEtAl2016} and made available online \citep{AhmeEtAl2015}.  As far as we know, problems in class SSLP were originally defined and solved in \citet{NtaiSen2005}. Section~\ref{sec:standardized_problem_families} briefly describes the problem families we focus on. From these families, we select instances that are reported as being hard to solve in the literature. We parameterize these instances so that we can generate samples of instances featuring similar characteristics. Section~\ref{sec:data_generation} details the data generation process and the input-output structures of our ML predictors. 

\subsection{Problem families}
\label{sec:standardized_problem_families}



The SSLP and SMKP specialize problem (P) described in Section~\ref{sec:intro} as follows:
\begin{enumerate}
\item 
In both SSLP and SMKP, $z$ is absent, i.e., $\mathcal{Z} = \emptyset$, and $\mathcal{Y}$ imposes only binary restrictions.

\item
In SSLP, $q_{\xi} \equiv q$, $W_{\xi} \equiv W$, $T_{\xi} \equiv T$, $\forall \xi$, i.e., all second-stage coefficients are deterministic, except the right-hand sides of some constraints.

\item In SMKP, $W_{\xi} \equiv W$, $h_{\xi} \equiv h$, $T_{\xi} \equiv T$, $\forall \xi$, i.e., all second-stage coefficients are deterministic, except those appearing in the objective.

\end{enumerate}

From SSLP \citep[Section~5.1][]{AnguEtAl2016}, we select problems SSLP(10, 50, 2000) and SSLP(15, 45, 15), where SSLP(a, b, c) features $a$ servers, $b$ clients and $c$ second-stage scenarios. According to the computations reported in \cite{AnguEtAl2016}, these are the most difficult to solve exactly among the problems in SSLP whose detailed statements have been made publicly available. The second stages of these problems are also among the most difficult to solve in SSLP. 

From SMKP \citep[Section~5.2][]{AnguEtAl2016}, we select problems indexed by 29 and 30, denoted SMKP(29) and SMKP(30). Similarly, these are the most difficult to solve exactly according to \cite{AnguEtAl2016} and the second stages of these problems are also among the most difficult in SMKP. We note that no solution was found to SMKP(30) in \cite{AnguEtAl2016} within the allocated time.

The difference in the nature of the problems belonging to SSLP and SMKP families is deep: in comparison with problems in SSLP, problems in SMKP feature considerably harder first stages and considerably easier integral second-stage problems. Clearly, the harder the second-stage problems are to solve exactly, the likelier is the expected gain in speed from ML-L-Shaped. Therefore, SSLP is a better candidate for our method than SMKP. In addition, because of the differences between these families, we expect the alternating cuts variant of ML-L-Shaped to perform better with SMKP than the standard integer L-shaped cut variant, whereas we expect the opposite to hold for SSLP. These opposing dominances are confirmed in our experiments.

Recall that our goal is to solve in short online computing time individual problems stemming from a distribution of problem instances sharing common characteristics. With this goal in mind, and the objective to draw statistical evidence from our numerical experiments, we turn our attention to the parameterization of these problems. We parameterize problems SSLP(10, 50, 2000) and SSLP(15, 45, 15) by allowing the individual deterministic capacities of the servers to vary. These capacities are specified in a deterministic subvector of $h_{\xi}$. Thus, whereas all servers of a problem have identical capacities in \cite{AnguEtAl2016}, these being respectively equal to 188 and 112 in SSLP(10, 50, 2000) and SSLP(15, 45, 15), we assume in effect that the problem instances that might be encountered are characterized by server-specific capacities ranging between 75 and 300. We shall let SSLPF(10, 50, 2000) and SSLPF(15, 45, 15) denote the derived families of problem instances obtained by parameterizing the original individual problems SSLP(10, 50, 2000) and SSLP(15, 45, 15).

Similarly, we parameterize problems SMKP(29) and SMKP(30) of \cite{AnguEtAl2016} by allowing the coefficients of the deterministic technology matrix $T$ and the deterministic right-hand side values $h$ appearing in the coupling constraints to vary, whereas the recourse matrix $W$ remains fixed across simulations. The individual coefficients of $T$ are integer-valued and range between 1 and 40 whereas the individual values in $h$ can vary freely, provided the resulting second-stage problems feature relatively complete recourse, given the values taken by $W$ and $T$. This gives rise to the two derived families of problem instances SMKPF(29) and SMKPF(30).

We also define two additional families: SSLPF(15, 45, 150) and SSLPF(15, 80, 15). These enable us to assess the effects of moderate increases in the complexity of second stage on the relative performances of the exact and approximate versions of the algorithms. These families feature exactly the same first-stage problem $(M)$ as that of SSLPF(15, 45, 15). Whereas SSLPF(15, 45, 150) shares the same recourse matrix $W$ with SSLPF(15, 45, 15) but features 150 instead of 15 scenarios in its second stage, SSLPF(15, 80, 15) shares the same 15 scenarios with SSLPF(15, 45, 15) but features 80 clients instead of 15. In SSLPF(15, 80, 15), the coefficients of the new recourse matrix $W$ are generated according to the same distribution as that originally used in \cite{AnguEtAl2016} for generating those of SSLP(15, 45, 15).

\subsection{Data generation}
\label{sec:data_generation}

Selection of data sets made up of (input, label) examples is critical for supervised ML. In the present context, an input is a statement of a second-stage problem (S), (RS) or (DRS) and a label is the  optimal solution thereof. The following details the choices governing the generation of our sets of examples.

\subsubsection{Families of SSLP instances}
\label{sec:families_of_SSLP_instances}

Generation of the data required for training and validating an ML predictor of the optimal objective value of second stage is performed in a similar fashion for each one of the standardized families SSLPF(10, 50, 2000), SSLPF(15, 45, 15), SSLPF(15, 45, 150) and SSLPF(15, 80, 15). Keeping all coefficients appearing in the second stage fixed except server capacities, the second-stage problems are simulated by pseudo-randomly sampling individual server capacities from independent discrete uniform distributions with support $[75, 300]$ and by pseudo-randomly sampling from independent discrete uniform distributions the values of the binary coupling variables shared by the first and second stages. Then, the corresponding optimal solution of second stage is computed exactly. Each such (problem statement, problem solution) pair constitutes a supervised example available to ML.

Problem statements pertaining to SSLPF(10, 50, 2000) are summarized by vectors in $\mathbb{N}^{20}$ (10 integral servers capacities + 10 coupling binaries), whereas problem statements pertaining to SSLPF(15, 45, 15), SSLPF(15, 45, 150) and SSLPF(15, 80, 15) are summarized by vectors in $\mathbb{N}^{30}$ (15 integral servers capacities + 15 coupling binaries). As expected, preliminary experiments have shown that the variant of Algorithm~\ref{alg:benders}with $\textit{isAlt}$ set to \emph{false} dominates the one where $\textit{isAlt}$ is set to \emph{true} in connection with the SSLP families. Hence, ML predictions $Q^{\textit{ML}}(x)$ required in building cuts (\ref{eq:integer_L_shaped_cut_heur}) suffice and all solutions pertaining to SSLP families can be described by points in $\mathbb{R}$.

Similarly to \cite{LarsEtAl2022a}, we deliberately generate very large numbers of examples to ensure that availability of data does not impinge on learning. For each one of the families, we produce a set of 1M examples. The data sets are partitioned according to proportions 64\%, 16\%, 20\% between training, validation and test sets. 

We also generate a second data set for the SSLPF(10, 50, 2000) family that we identify as SSLPF-indx(10, 50, 2000). Examples in this set comprise the same problem statements as those of SSLPF(10, 50, 2000). However, instead of calculating the solution as an expectation over all second-stage scenarios, we calculate it for a single randomly selected scenario. It is considerably faster to generate each example in this fashion: the time required for producing 1M examples of the latter is approximately 1/80th of that required to generate 1M examples of SSLPF(10, 50, 2000). \citet{LarsEtAl2022a} have argued and provided evidence that the resulting ML predictors should have predictive properties close to those of ML predictors computed from examples where the solutions reflect all scenarios, calling this phenomenon ``implicit aggregation'' and the resulting predictor ``implicit predictor''. We later compare the performances achieved in approximating first-stage solutions through ML predictors originating from SSLPF(10, 50, 2000) and SSLPF-indx(10, 50, 2000) data.

Average generation time per example of each SSLP family is below 0.6 second for all families except SSLPF(10, 50, 2000) where the average is 3.24 seconds. We report additional detail in Appendix~A (Table~\ref{tab:average_generation_time_per_example}).

\subsubsection{Families of SMKP instances}
\label{families_of_smkp_instances}

Generation is performed in a distinct manner for the SMKP families. If the statements of the problems in SMKPF(29) and SMKPF(30) were to be described by stating separately the values of each binary coupling variable, each coefficient of the technology matrix $T$ and each element of the right-hand side $h$, then a vector in $\mathbb{N}^{600}$ -- since there are 100 individual coupling variables in $x$ and $T$ is of dimension $(5\times100)$ -- and a vector in $\mathbb{R}^{5}$ -- since $h$ is of dimension $(5\times1)$ -- would be required. However, problems in SMKPF(29) and SMKPF(30) can be stated compactly with a vector in $\mathbb{R}^{5}$ since each problem is fully described by the real vector $h + Tx$ of dimension $(5\times1)$. Reducing the input in this fashion makes it far less demanding for ML. In general, similar input reductions are applicable in contexts where a detailed problem description stating all individual input elements would hamper the learning process.

As expected, our experiments show that, in contrast with families in SSLP, ML-L-Shaped with $\textit{isAlt}$ set to \emph{true} should be used with the SMKP families.  Therefore, we train and validate two ML predictors. One yields the predictor of the expected objective value of the original second-stage subproblem (S), $Q^{\textit{ML}}(x)$, required in building cuts (\ref{eq:integer_L_shaped_cut_heur}). Another outputs (i) $\widetilde{Q}^{\textit{ML}}(x^*)$, (ii) $\lbrace \mathbb{E}_{\xi} [\phi_{\xi}]\rbrace ^{ML}$ and (iii) $\lbrace \mathbb{E}_{\xi} [\mathbf{1}' \psi_{\xi}]\rbrace ^{ML}$. (i), (ii) and (iii) respectively predict the optimal expected objective value of the continuously-relaxed second-stage subproblem (RS), the expectation of the dual variables attached to the coupling constraints and the expectation of the summation of the dual variables attached to the upper bounds (equal to 1) of all second-stage variables. In regard to SMKPF(29) and SMKPF(30), (i) has value in $\mathbb{R}$, (ii) has value in $\mathbb{R}^{5}$ and (iii) has value in $\mathbb{R}$. Hence, the output can be completely described by a vector in $\mathbb{R}^{7}$. In contrast, an extensive description of the solution to the optimal solution of the continuously-relaxed second stage that would account for every dual value and every scenario would involve several hundred elements. Summarizing the statement of the solution in this manner is far less demanding for ML. Similar output reductions are applicable in contexts where a description of all individual output elements would be overly detailed. 

For each of SMKPF(29) and SMKPF(30), two data sets are generated. The two data sets share the same reduced problem statements. However, one data set comprises the solutions to the integral second-stage problem (S), whereas the other comprises the reduced descriptions of solutions to its continuous relaxation (DRS). All four resulting data sets are partitioned according to proportions 64\%, 16\%, 20\% between training, validation and test sets. The average generation time per example is less than 0.13 second. We report additional detail in Table~\ref{tab:average_generation_time_per_example} of Appendix~A.

\section{ML predictors} \label{sec:ml_approximation}

Based on the data sets generated according to Section~\ref{sec:data_generation}, we build the ML-predictors yielding the approximate second-stage solutions called at Steps~\ref{a:contPred}~and \ref{heur_is_not_oppor2} of Algorithm~\ref{alg:heuristic_callback_v2}. We proceed to detail their characteristics, the process followed in their training/validation and to report their predictive performance over their respective test sets.

As our primary goal is to assess the usefulness of the ML-L-Shaped algorithm in solving problems of class (P), we implement standard feed-forward neural network architectures. These generic ML approximators feature a well-documented performance (in terms of speed and accuracy) and are known to scale well to wide ranges of input and output lengths and values. Similarly, our hyperparameter search, inclusive of the characteristics of the feed-forward networks, remains limited in view of our goal of illustrating the implementation of our approximation methods and probing their potential benefits rather than extracting maximal predictive accuracy or data-efficiency.
We employ GPUs for training/validation and generating predictions. Whereas these operations can run on CPUs, the gains in performance offered by GPUs make their use practically indispensable, despite the additional challenges associated, for example, with the calls from CPU threads occurring at prediction time at Steps~\ref{a:contPred}~and~\ref{heur_is_not_oppor2} of Algorithm~\ref{alg:heuristic_callback_v2}.

In view of their heterogeneous ranges, the non-binary inputs of the predictors for the SSLP families are rescaled in $[0, 1]$ to ensure better numerical conditioning during training. In contrast, the transformation applied to the inputs of the predictors for the SMKP families makes rescaling unnecessary. 

We minimize L1 error over the training set with stochastic mini-batch gradient descent equipped with Adam learning rate adaptation and mini-batch size equal to 128. Weighting inversely proportional to the sample averages of the output values measured on the training set is applied to the individual L1 errors of the networks outputting multiple values when calculating the training and validation errors. 

Early stopping is applied based on L1 error (weighted when there are multiple outputs) over the full validation set with a forgiving patience of at least several hundred epochs. Networks outputting a single value $Q^{\textit{ML}}(x)$ (SSLP and SMKP families) are equipped with 10 hidden layers of 800 units each. Networks outputting multiple values (SMKP families) are equipped with 15 hidden layers of 1,000 units each. All units except those in last hidden and output layers are fitted with rectified linear activations. Units in last hidden and output layers are fitted with linear activations. ML training and validation is performed with Python PyTorch 1.71 on a single Nvidia RTX 3090 GPU.

\begin{table}[htbp]
  \centering
  \resizebox{\columnwidth}{!}{
    \begin{tabular}{||l||c|c|c|c|c|c||}
    \hline
    \textbf{Problem family} &  IP/LP &  Input &  \# Hid. &  units/ &  Output &  Abs. rel.  \\
	 &   &  length &  layers &  hid. layer &  length &  error [\%] \\
    \hline
    \hline
    \textbf{SSLPF(10,50,2000)} & IP & 20 & 10 & 800 & 1 & 0.87  \\
    \hline
    \textbf{SSLPF-indx(10,50,2000)} & IP  & 20 & 10 & 800 & 1 & 5.31  \\
    \hline
    \textbf{SSLPF(15,45,15)} & IP & 30 & 10 & 800 & 1 & 0.23  \\
    \hline
    \textbf{SSLPF(15,45,150)} & IP & 30 & 10 & 800 & 1 & 0.12  \\
    \hline
    \textbf{SSLPF(15,80,15)} & IP & 30 & 10 & 800 & 1 & 0.40  \\
    \hline
    \hline
    \textbf{SMKPF(29)} & IP &  5 & 10 & 800 & 1 & 0.071  \\
    \hline
    \textbf{SMKPF(29)} & LP &  5 & 15 & 1000 & 7 & 6.64  \\
    \hline
    \textbf{SMKPF(30)} & IP &  5 & 10 & 800 & 1 & 0.072  \\
    \hline
    \textbf{SMKPF(30)} & LP &  5 & 15 & 1000 & 7 & 7.41  \\
    \hline
     \multicolumn{7}{l}{\small IP, LP: output is solution of integral or relaxed 2nd stage problem.} \\
     \multicolumn{7}{l}{\small Abs. rel. error: average absolute relative prediction error made on ML test set.} \\
     \multicolumn{7}{l}{The test set is same for SSLPF-indx(10,50,2000) and SSLPF(10,50,2000).} \\
    \end{tabular}
    }
    \caption{Performance of ML predictors}
    \label{tab:ml_predictors}
\end{table}

Table~\ref{tab:ml_predictors} reports the characteristics of the ML predictors and their performance on their respective test sets. Total of training and validation times is reported for each predictor in Table~\ref{tab:training_and_validation_times} of Appendix~A. Crucially, the computational advantage presented by the ML-L-Shaped method hinges on the time within which a prediction can be generated. For every feed-forward network considered here, this time is nearly constant and at most equal to a few milliseconds. The predictive performance appears to be good or excellent. Measured in absolute relative percentage error, it is better than 1\% for models that output $Q^{\textit{ML}}(x)$ and are trained on the most accurate data (see rows with IP label). A larger error (5.31\% vs 0.87\%) ensues when performing implicit aggregation \citep{LarsEtAl2022a}. The tasks related to the continuous cuts (\ref{eq:continuous_L_shaped_cut_heur}) are harder (see SMKP rows with LP label), and errors vary in this context between 6\% and 8\%.

\section{Experimental results}
\label{sec:detailed_analysis}

This section compares the performances achieved by the ML-L-Shaped heuristic and the best performing exact method. The results we report are based on independently generated instances of the SSLP and SMKP families. Section~\ref{sec:comparison_with_alternative_heuristic} is dedicated to a comparison between ML-L-Shaped and an alternative heuristic (PH).

Our comparisons pertain to computing time, optimality gap, number of first-stage nodes, number and duration of integral second-stage problems, number and duration of continuously-relaxed second-stage problems. We also present evidence regarding the use of heuristic solutions as warm-start incumbents. All computations reported in this section are performed on an Intel i9-10980XE processor using the Java programming language adjoined with the Java version of CPLEX, version 12.10. ML predictions are called through the Java bindings of TorchScript and performed on a single Nvidia RTX 3090 GPU.

Tables~\ref{tab:computation_times}~to~\ref{tab:total_computation_times_with_feedback_from_approximate_solution} detail the results of our experiments for both SSLP and SMKP families whereas the relevant analyses are split between Sections~\ref{sec:families_of_sslp_instances}~and~\ref{sec:families_of_smkp_instances}. We open with a number of background remarks.
\emph{First}, the shorthand notations \emph{Std-L} and \emph{Alt-L} stand for respectively the exact standard integer L-shaped method and the one with alternating cuts. Similarly, \emph{ML-Std-L} and \emph{ML-Alt-L} designate the two variants of the ML-L-Shaped method.
\emph{Second}, for each one of the SSLP and SMPK families, the reported statistics are calculated from 100 independently sampled instances that have not been previously involved either in ML or in preliminary experiments or calibrations. This guards against a well-known source of overoptimism. \emph{Third}, our results about SSLPF(10, 50, 2000), SSLPF(15, 45, 15), SMKPF(29), SMKPF(30) are aligned with those reported in \cite{AnguEtAl2016}: We find out that Alt-L achieves computing times that are smaller or equal to those of Std-L. Hence, Alt-L is reported in the tables as the benchmark for any comparison with ML-L-Shaped. Moreover, the hierarchy of our computing times between SSLPF(10, 50, 2000), SSLPF(15, 45, 15), SMKPF(29), SMKPF(30) is similar to that reported in \cite{AnguEtAl2016}. \emph{Fourth}, the shift coefficients are as follows: $\mu = 1.0$ in connection with instances belonging to the SSLP families, $\mu = 0.98$ and $\nu = 0.95$ in connection with instances belonging to the SMKP families. We encountered one instance of SMKPF(29) and one instance of SMKPF(30) where these values were too large to yield a feasible solution. For each of these two instances, we reapplied ML-L-Shaped with decreasing values of $\mu$ and $\nu$, as indicated in Section~\ref{sec:decomposition_algorithms}, reaching solutions with $\mu = \nu = 0.7$.

\subsection{Families of SSLP instances}
\label{sec:families_of_sslp_instances}
This section focuses on the results related to the SSLP families. Those are reported in the first five rows of Tables~\ref{tab:computation_times}~to~\ref{tab:total_computation_times_with_feedback_from_approximate_solution} and reflect the application of the ML-Std-L variant of ML-L-Shaped. That ML-Std-L dominates ML-Alt-L on SSLP instances was established by independent, preliminary experiments, as expected.

\paragraph{Computing times.}
\label{par:sslp_computation_times}
As evidenced in Table~\ref{tab:computation_times}, ML-Std-L achieves computing times far smaller than those required by Alt-L. For instance, the average ratios of computing times over 100 generated instances of families SSLPF(10, 50, 2000) and SSLPF(15, 45, 15) are respectively equal to 0.57\% and 8.67\% (representing speed-ups of 175x and 12x). Noticing that the second-stage problems in family SSLPF(15, 45, 15) are quite simple, leading to an average overall exact solution time of 5.25~seconds, we present complementary results for the families SSLPF(15, 45, 150) and SSLPF(15, 80, 15) whose second stages are moderately more complex. The corresponding average ratios of computing times, respectively equal to 1.57\% and 7.30\% (representing speed-ups of 64x and 14x), illustrate how ML-Std-L becomes more advantageous as the complexity of the second-stage problems increases. As expected, since ML prediction times are nearly invariant, the ratio of computing times obtained for SSLPF-indx(10, 50, 2000), equal to 0.52\%, is nearly equal to that of SSLPF(10, 50, 2000).

\begin{table}[htbp]
  \centering
  \resizebox{\columnwidth}{!}{
    \begin{tabular}{||l||rrr|r||rrr|r||rrr|r||}
    \hline
    \multirow{3}{*}{\textbf{Problem family}}
    & \multicolumn{4}{c||}{Alt-L} & \multicolumn{4}{c||}{ML-L-Shaped} & \multicolumn{4}{c||}{ML-L-Shaped/Alt-L ratio}    \\
	& \multicolumn{3}{c|}{Quantiles} &  &  \multicolumn{3}{c|}{Quantiles}  &   & \multicolumn{3}{c|}{Quantiles} &    \\
	&  0.05 & 0.5 & 0.95 & Avg & 0.05 & 0.5 & 0.95 & Avg & 0.05 & 0.5 & 0.95 & Avg \\
    \hline
    \hline
    \textbf{SSLPF(10,50,2000)} & 131.63 & 150.73 & 186.04 & 156.06& 0.52& 0.64& 2.59& 0.93 & 0.34\% & 0.43\% & 1.49\% & 0.57\% \\   
    &       &       &       & (2.89)  &       &       &       & (0.08)  &       &       &       & (0.03\%) \\   
    \hline
    \textbf{SSLPF-indx(10,50,2000)} & 131.63 & 150.73 & 186.04 & 156.06 & 0.46  & 0.60  & 2.39  & 0.85  & 0.31\% & 0.42\% & 1.42\% & 0.52\% \\    
    &       &       &       & (2.89)  &       &       &       & (0.08)  &       &       &       & (0.03\%) \\  
    \hline
    \textbf{SSLPF(15,45,15)}       & 4.28  & 5.00  & 7.12  & 5.25  & 0.37  & 0.43  & 0.60  & 0.45  & 7.55\% & 8.68\% & 9.99\% & 8.67\% \\    
    &       &       &       & (0.11)  &       &       &       & (0.01)  &       &       &        & (0.07\%) \\  
    \hline
    \textbf{SSLPF(15,45,150)}        & 27.96 & 34.16 & 43.93 & 34.96 & 0.39  & 0.53  & 0.79  & 0.55  & 1.32\% & 1.51\% & 2.13\% & 1.57\% \\    
    &       &       &       & (0.57)  &       &       &       & (0.01)  &       &       &        & (0.03\%) \\    
    \hline
    \textbf{SSLPF(15,80,15)}       & 36.20 & 47.24 & 130.70 & 58.55 & 2.32  & 3.22  & 13.90 & 4.12  & 3.60\% & 6.67\% & 14.54\% & 7.30\%  \\  
    &       &       &       & (3.15)  &       &       &       & (0.29)  &       &       &       & (0.36\%) \\    
    \hline   \hline
    \textbf{SMKPF(29)} &       151.82 & 742.89 & 3996.73 & 1233.74 & 28.02 & 109.12 & 510.91 & 174.55 & 4.87\% & 16.51\% & 51.22\% & 19.99\% \\    
    &       &       &       & (123.09) &       &       &       & (18.05) &       &       &       & (1.47\%) \\   
    \hline
    \textbf{SMKPF(30)} &      212.54 & 1459.19 & 6785.80 & 2132.51 & 37.91 & 194.34 & 935.81 & 329.82 & 3.66\% & 15.70\% & 57.38\% & 19.34\% \\    
    &       &       &       & (254.70) &       &       &       & (55.97) &       &       &       & (1.46\%) \\
    \hline
    \multicolumn{13}{l}{\small Standard error of estimate is reported between parentheses.} \\
    \end{tabular}    
    }
    \caption{Computing times (seconds)}
    \label{tab:computation_times}
\end{table}

\paragraph{First-stage objective values.}
\label{par:sslp_first_stage_objective_values}

As evidenced by Table~\ref{tab:first_stage_values}, the performance of ML-Std-L in regard to the values achieved for the objective of first stage is excellent. The average optimality gaps over the instances of families SSLPF(10, 50, 2000), SSLPF(15, 45, 15), SSLPF(15, 45, 150) and SSLPF(15, 80, 15), vary between nearly zero and 1.94\%. In addition, equal to 2.61\%, the small first-stage optimality gap yielded by the implicit ML predictor of the value of the second-stage problems of SSLPF(10, 50, 2000) (see SSLPF-indx(10, 50, 2000)) provides evidence in favor of implicit second-stage predictors when generating second-stage data based on all scenarios is highly time-consuming. (Recall that generating SSLPF-indx(10, 50, 2000) requires the solution to the second-stage problem of a single scenario.)

\begin{table}[htbp]
  \centering
  \resizebox{\columnwidth}{!}{
    \begin{tabular}{||l||rrr|r||rrr|r||rrr|r||}
    \hline
    \multirow{3}{*}{\textbf{Problem family}}
    & \multicolumn{4}{c||}{Alt-L} & \multicolumn{4}{c||}{ML-L-Shaped} & \multicolumn{4}{c||}{Optimality gap}    \\
	& \multicolumn{3}{c|}{Quantiles} &  &  \multicolumn{3}{c|}{Quantiles}  &   & \multicolumn{3}{c|}{Quantiles} &    \\
	&  0.05 & 0.5 & 0.95 & Avg & 0.05 & 0.5 & 0.95 & Avg & 0.05 & 0.5 & 0.95 & Avg \\
	\hline
    \hline
    \textbf{SSLPF(10,50,2000)} & -353.85 & -347.14 & -333.17 & -345.60 & -353.85 & -347.14 & -333.17 & -345.58 & 0.000\% & 0.000\% & 0.036\% & 0.006\% \\   
    &       &       &       & (0.70)  &       &       &       & (0.70)  &       &       &       & (0.003\%) \\  
    \hline
    \textbf{SSLPF-indx(10,50,2000)} & -353.85 & -347.14 & -333.17 & -345.60 & -348.92 & -339.00 & -315.25 & -336.57 & 0.000\% & 2.173\% & 7.850\% & 2.609\% \\    
    &       &       &       & (0.70)  &       &       &       & (1.06)  &       &       &       & (0.242\%) \\
    \hline
    \textbf{SSLPF(15,45,15)}       & -313.20 & -308.74 & -296.62 & -307.15 & -313.20 & -308.67 & -296.48 & -306.95 & 0.000\% & 0.000\% & 0.608\% & 0.064\% \\    
    &       &       &       & (0.59)  &       &       &       & (0.60)  &       &       &       & (0.019\%) \\
    \hline
    \textbf{SSLPF(15,45,150)}        & -314.16 & -306.42 & -294.74 & -305.89 & -314.15 & -306.27 & -281.51 & -300.01 & 0.000\% & 0.000\% & 1.021\% & 1.943\% \\    
    &       &       &       & (0.66)  &       &       &       & (3.61)  &       &       &       & (1.150\%) \\   
    \hline
    \textbf{SSLPF(15,80,15)}       & -632.21 & -614.67 & -584.42 & -613.43 & -630.52 & -614.67 & -584.42 & -612.97 & 0.000\% & 0.000\% & 0.538\% & 0.075\%   \\
    &       &       &       & (1.34)  &       &       &       & (1.34)  &       &       &       & (0.018\%) \\
    \hline   \hline
    \textbf{SMKPF(29)} &       7736.92 & 8176.1 & 8567.36 & 8155.42 & 7736.92 & 8178.5 & 8570.25 & 8156.18 & 0.000\% & 0.000\% & 0.050\% & 0.009\% \\
    &       &       &       & (27.60) &       &       &       & (25.58) &       &       &       & (0.002\%) \\ 
    \hline
    \textbf{SMKPF(30)} &      8288.24 & 8790.53 & 9160.83 & 8754.33 & 8228.24 & 8790.53 & 9164.16 & 8754.73 & 0.000\% & 0.000\% & 0.027\% & 0.005\% \\
    &       &       &       & (28.18) &       &       &       & (28.19) &       &       &       & (0.001\%) \\
    \hline
    \multicolumn{13}{l}{\small Standard error of estimate is reported between parentheses.} \\
    \end{tabular}    
    }
    \caption{First-stage values and optimality gaps}
    \label{tab:first_stage_values}
\end{table}

\paragraph{Number of first-stage problems.}
\label{par:sslp_number_of_first_stage_problems}

Table~\ref{tab:number_of_first_stage_nodes} indicates that the average ratio of the numbers of first-stage nodes between approximate and exact solutions varies over the SSLP families between 50\%  and 76\%, thus excluding simplification of the first-stage problems as the main source of reduction in computing times between Alt-L and ML-Std-L.

\begin{table}[htbp]
  \centering
  \resizebox{\columnwidth}{!}{
    \begin{tabular}{||l||rrr|r||rrr|r||rrr|r||}
    \hline
    \multirow{3}{*}{\textbf{Problem family}}
    & \multicolumn{4}{c||}{Alt-L} & \multicolumn{4}{c||}{ML-L-Shaped} & \multicolumn{4}{c||}{ML-L-Shaped/Alt-L ratio}    \\
	& \multicolumn{3}{c|}{Quantiles} &  &  \multicolumn{3}{c|}{Quantiles}  &   & \multicolumn{3}{c|}{Quantiles} &    \\
	&  0.05 & 0.5 & 0.95 & Avg & 0.05 & 0.5 & 0.95 & Avg & 0.05 & 0.5 & 0.95 & Avg \\
	\hline
    \hline
    \textbf{SSLPF(10,50,2000)} & 5.51E+02 & 6.21E+02 & 7.16E+02 & 6.26E+02 & 3.85E+02 & 4.27E+02 & 4.83E+02 & 4.29E+02 & 64.42\% & 68.50\% & 72.60\% & 68.52\%  \\
    &       &       &       & (4.96E+00) &       &       &       & (3.39E+00) &       &       &       & (0.24\%) \\
    \hline
    \textbf{SSLPF-indx(10,50,2000)} & 5.51E+02 & 6.21E+02 & 7.16E+02 & 6.26E+02 & 3.29E+02 & 4.32E+02 & 5.14E+02 & 4.27E+02 & 56.93\% & 68.97\% & 78.33\% & 68.26\% \\
    &       &       &       & (4.96E+00) &       &       &       & (4.92E+00) &       &       &       & (0.67\%) \\
    \hline
    \textbf{SSLPF(15,45,15)}       & 1.82E+03 & 2.05E+03 & 2.64E+03 & 2.13E+03 & 1.33E+03 & 1.58E+03 & 2.03E+03 & 1.62E+03 & 68.82\% & 76.20\% & 82.52\% & 75.95\% \\
    &       &       &       & (2.62E+01) &       &       &       & (2.30E+01) &       &       &       & (0.41\%) \\
    \hline
    \textbf{SSLPF(15,45,150)}        & 1.83E+03 & 2.15E+03 & 2.79E+03 & 2.19E+03 & 1.16E+03 & 1.61E+03 & 2.04E+03 & 1.59E+03 & 64.28\% & 74.52\% & 82.02\% & 72.97\% \\
    &       &       &       & (2.88E+01) &       &       &       & (3.33E+01) &       &       &       & (1.13\%) \\  
    \hline
    \textbf{SSLPF(15,80,15)}    & 1.27E+04 & 1.42E+04 & 1.73E+04 & 1.44E+04 & 5.56E+03 & 6.76E+03 & 1.45E+04 & 7.29E+03 & 40.75\% & 48.12\% & 81.04\% & 49.96\% \\
    &       &       &       & (1.32E+02) &       &       &       & (2.20E+02) &       &       &       & (1.02\%) \\
    \hline   \hline
    \textbf{SMKPF(29)} &       1.20E+06 & 6.04E+06 & 3.12E+07 & 9.59E+06 & 2.21E+05 & 8.21E+05 & 3.91E+06 & 1.35E+06 & 4.47\% & 15.42\% & 49.92\% & 19.16\% \\
    &       &       &       & (9.18E+05) &       &       &       & (1.43E+05) &       &       &       & (1.40\%) \\
    \hline
    \textbf{SMKPF(30)} &      1.71E+06 & 1.17E+07 & 4.47E+07 & 1.59E+07 & 2.72E+05 & 1.54E+06 & 7.39E+06 & 2.51E+06 & 3.65\% & 14.53\% & 56.20\% & 18.62\% \\
    &       &       &       & (1.49E+06) &       &       &       & (3.90E+05) &       &       &       & (1.44\%) \\
    \hline
    \multicolumn{13}{l}{\small Standard error of estimate is reported between parentheses.} \\
    \end{tabular}    
    }
    \caption{Number of first-stage nodes}
    \label{tab:number_of_first_stage_nodes}
\end{table}

\paragraph{Second-stage integral and relaxed problems.}
\label{par:sslp_second_stage_integral_and_relaxed_problems}

Tables~\ref{tab:integral_second_stage_problems}~and~\ref{tab:relaxed_second_stage_problems} report the numbers of integral and relaxed second-stage problems and the total times spent in the latter. Their content reflects the parallelism taking place in CPLEX's branch-and-cut computations, which entails that the net computing times appearing in Table~\ref{tab:computation_times} cannot be directly inferred from it. Average times per individual integral problem can, however, be estimated by dividing the average total time and the average number of problems (reported in lower and upper panels of Table~\ref{tab:integral_second_stage_problems}). Average times per individual relaxed problem can be estimated in a similar fashion from Table~\ref{tab:relaxed_second_stage_problems}.

By construction, Std-L only requires the solution of integral second-stage problems. Table~\ref{tab:integral_second_stage_problems} shows that, when it is applied to the SSLP families, ML-Std-L must compute at most a few thousand approximate integer L-shaped cuts and that the latter are generated within a few milliseconds on average, including the time required to predict the value of second stage with ML. In contrast, Table~\ref{tab:relaxed_second_stage_problems} shows that Alt-L must solve exactly a comparable number of relaxed second-stage problems at a cost of several milliseconds each. In addition, 
according to Table~\ref{tab:integral_second_stage_problems}, it must solve exactly several integral second-stage problems at a cost of one or more seconds each.

The latter findings point to the main source of the reduction in computing time achieved by ML-Std-L in comparison to the fastest exact algorithm available (Alt-L), when the number of first-stage nodes is in the order of thousands as in the SSLP families: The high speed with which predictions of the value of the integral second-stage problem are generated with ML is sufficient to offset a larger number of integral second-stage problems in comparison with Alt-L. This speed is even higher than the speed at which the exact solutions of the relaxed second-stage problems can be computed by Alt-L.

\begin{table}[htbp]
  \centering
  \resizebox{\columnwidth}{!}{
    \begin{tabular}{||l||rrr|r||rrr|r||rrr|r||}
    \hline
    \multirow{3}{*}{\textbf{Problem family}}
    & \multicolumn{4}{c||}{Alt-L} & \multicolumn{4}{c||}{ML-L-Shaped} & \multicolumn{4}{c||}{ML-L-Shaped/Alt-L ratio}    \\
	& \multicolumn{3}{c|}{Quantiles} &  &  \multicolumn{3}{c|}{Quantiles}  &   & \multicolumn{3}{c|}{Quantiles} &    \\
	&  0.05 & 0.5 & 0.95 & Avg & 0.05 & 0.5 & 0.95 & Avg & 0.05 & 0.5 & 0.95 & Avg \\
	\hline
    \hline
    \multicolumn{13}{||c||}{} \\
    \multicolumn{13}{||c||}{\textbf{Number of problems}} \\
    \hline
	\textbf{SSLPF(10,50,2000)} & 56.00 & 65.00 & 82.90 & 66.73 & 381.00 & 417.00 & 455.95 & 419.17 & 514.06\% & 643.20\% & 756.01\% & 636.07\% \\
    &       &       &       & (0.88)  &       &       &       & (2.73)  &       &       &       & (7.14\%) \\
    \hline
    \textbf{SSLPF-indx(10,50,2000)} & 56.00 & 65.00 & 82.90 & 66.73 & 327.10 & 402.00 & 467.00 & 401.16 & 475.16\% & 621.40\% & 733.87\% & 608.65\% \\
    &       &       &       & (0.88)  &       &       &       & (3.98)  &       &       &       & (8.06\%) \\
    \hline
    \textbf{SSLPF(15,45,15)}       & 52.00 & 63.00 & 79.95 & 64.06 & 1075.05 & 1243.50 & 1601.55 & 1289.88 & 1563.66\% & 1973.50\% & 2749.83\% & 2052.07\% \\
    &       &       &       & (0.96)  &       &       &       & (17.59) &       &       &       & (38.41\%) \\
    \hline
    \textbf{SSLPF(15,45,150)}        & 57.05 & 70.00 & 86.00 & 70.29 & 1006.25 & 1259.00 & 1582.85 & 1266.55 & 1289.20\% & 1791.57\% & 2556.04\% & 1828.42\% \\
    &       &       &       & (0.84)  &       &       &       & (24.91) &       &       &       & (42.47\%) \\
    \hline
    \textbf{SSLPF(15,80,15)}    & 37.10 & 72.00 & 89.95 & 70.26 & 4952.95 & 6052.00 & 13922.15 & 6641.81 & 6587.86\% & 8504.23\% & 18138.67\% & 9875.65\% \\
    &       &       &       & (1.48)  &       &       &       & (221.40) &       &       &       & (362.88\%) \\
    \hline   \hline
    \textbf{SMKPF(29)} &       35.00 & 64.00 & 177.65 & 77.12 & 13.00 & 20.00 & 33.00 & 21.54 & 13.09\% & 31.23\% & 62.96\% & 33.66\% \\
    &       &       &       & (4.47)  &       &       &       & (0.74)  &       &       &       & (1.51\%) \\
    \hline
    \textbf{SMKPF(30)} &      36.15 & 71.00 & 238.65 & 99.65 & 14.00 & 24.50 & 42.90 & 29.66 & 9.78\% & 31.55\% & 77.69\% & 35.43\% \\
    &       &       &       & (8.88)  &       &       &       & (4.76)  &       &       &       & (2.10\%) \\
    \hline
    \multicolumn{13}{||c||}{} \\
    \multicolumn{13}{||c||}{\textbf{Total time spent (milliseconds)}} \\
    \hline
    \textbf{SSLPF(10,50,2000)} & 26439.20 & 36704.00 & 63863.85 & 42256.11 & 568.05 & 622.00 & 709.55 & 632.82 & 1.02\% & 1.77\% & 2.32\% & 1.69\% \\
    &       &       &       & (2867.43) &       &       &       & (6.44)  &       &       &       & (0.04\%) \\
    \hline
    \textbf{SSLPF-indx(10,50,2000)} & 26439.20 & 36704.00 & 63863.85 & 42256.11 & 530.25 & 643.50 & 772.35 & 647.42 & 1.09\% & 1.73\% & 2.34\% & 1.73\% \\
    &       &       &       & (2867.43) &       &       &       & (7.89)  &       &       &       & (0.04\%) \\
    \hline
    \textbf{SSLPF(15,45,15)}       & 269.50 & 511.50 & 1434.75 & 633.63 & 1606.10 & 1845.50 & 2268.75 & 1876.13 & 141.95\% & 357.32\% & 684.90\% & 370.39\% \\
    &       &       &       & (34.56) &       &       &       & (21.33) &       &       &       & (16.43\%) \\
    \hline
    \textbf{SSLPF(15,45,150)}        & 2140.10 & 3851.00 & 7380.75 & 4097.17 & 1539.75 & 1891.50 & 2339.10 & 1886.55 & 27.87\% & 48.97\% & 87.26\% & 51.92\% \\
    &       &       &       & (166.85) &       &       &       & (33.91) &       &       &       & (1.82\%) \\
    \hline
    \textbf{SSLPF(15,80,15)}    & 2492.65 & 9765.00 & 84545.40 & 20430.97 & 4882.90 & 5599.50 & 10086.90 & 5921.83 & 8.76\% & 57.08\% & 245.03\% & 90.41\% \\
    &       &       &       & (2999.29) &       &       &       & (131.59) &       &       &       & (10.39\%) \\
    \hline   \hline
    \textbf{SMKPF(29)} &       323.20 & 1390.50 & 13930.55 & 3247.54 & 13.05 & 21.00 & 54.70 & 25.02 & 0.19\% & 1.53\% & 7.94\% & 2.33\% \\
    &       &       &       & (483.22) &       &       &       & (1.26)  &       &       &       & (0.22\%) \\
    \hline
    \textbf{SMKPF(30)} &      367.10 & 1964.50 & 12432.40 & 4633.48 & 13.00 & 27.00 & 62.95 & 35.08 & 0.28\% & 1.21\% & 7.37\% & 2.30\% \\
    &       &       &       & (1269.85) &       &       &       & (5.67)  &       &       &       & (0.23\%) \\
    \hline
    \multicolumn{13}{l}{\small Standard error of estimate is reported between parentheses.} \\
    \end{tabular}    
    }
    \caption{Integral second-stage problems}
    \label{tab:integral_second_stage_problems}
\end{table}

\begin{table}[htbp]
  \centering
  \resizebox{\columnwidth}{!}{
    \begin{tabular}{||l||rrr|r||rrr|r||rrr|r||}
    \hline
    \multirow{3}{*}{\textbf{Problem family}}
    & \multicolumn{4}{c||}{Alt-L} & \multicolumn{4}{c||}{ML-L-Shaped} & \multicolumn{4}{c||}{ML-L-Shaped/Alt-L ratio}    \\
	& \multicolumn{3}{c|}{Quantiles} &  &  \multicolumn{3}{c|}{Quantiles}  &   & \multicolumn{3}{c|}{Quantiles} &    \\
	&  0.05 & 0.5 & 0.95 & Avg & 0.05 & 0.5 & 0.95 & Avg & 0.05 & 0.5 & 0.95 & Avg \\
	\hline
    \hline
    \multicolumn{13}{||c||}{} \\
    \multicolumn{13}{||c||}{\textbf{Number of problems}} \\
    \hline
    \textbf{SSLPF(10,50,2000)} & 364.00 & 406.00 & 460.95 & 408.13 & 0.00  & 0.00  & 0.00  & 0.00  & 0.00\% & 0.00\% & 0.00\% & 0.00\% \\
    &       &       &       & (2.97)  &       &       &       & (0.00)  &       &       &       & (0.00\%) \\
    \hline
    \textbf{SSLPF-indx(10,50,2000)} & 364.00 & 406.00 & 460.95 & 408.13 & 0.00  & 0.00  & 0.00  & 0.00  & 0.00\% & 0.00\% & 0.00\% & 0.00\% \\
    &       &       &       & (2.97)  &       &       &       & (0.00)  &       &       &       & (0.00\%) \\
    \hline
    \textbf{SSLPF(15,45,15)}       & 933.15 & 1047.50 & 1347.20 & 1084.63 & 0.00  & 0.00  & 0.00  & 0.00  & 0.00\% & 0.00\% & 0.00\% & 0.00\% \\
    &       &       &       & (14.04) &       &       &       & (0.00)  &       &       &       & (0.00\%) \\
    \hline
    \textbf{SSLPF(15,45,150)}        & 924.70 & 1089.00 & 1392.25 & 1115.20 & 0.00  & 0.00  & 0.00  & 0.00  & 0.00\% & 0.00\% & 0.00\% & 0.00\% \\
    &       &       &       & (15.40) &       &       &       & (0.00)  &       &       &       & (0.00\%) \\
    \hline
    \textbf{SSLPF(15,80,15)}    & 5432.45 & 6234.00 & 7711.50 & 6308.44 & 0.00  & 0.00  & 0.00  & 0.00  & 0.00\% & 0.00\% & 0.00\% & 0.00\% \\
    &       &       &       & (65.66) &       &       &       & (0.00)  &       &       &       & (0.00\%) \\
    \hline   \hline
    \textbf{SMKPF(29)} &       55.15 & 108.50 & 330.20 & 137.04 & 17.00 & 24.00 & 48.55 & 29.21 & 8.11\% & 22.88\% & 51.11\% & 27.95\% \\
    &       &       &       & (8.79)  &       &       &       & (3.60)  &       &       &       & (4.17\%) \\
    \hline
    \textbf{SMKPF(30)} &      55.05 & 124.00 & 420.80 & 174.40 & 18.00 & 29.00 & 48.00 & 34.62 & 6.68\% & 21.56\% & 59.62\% & 25.59\% \\
    &       &       &       & (16.52) &       &       &       & (4.78)  &       &       &       & (1.64\%) \\
    \hline
    \multicolumn{13}{||c||}{} \\
    \multicolumn{13}{||c||}{\textbf{Total time spent (milliseconds)}} \\
    \hline
    \textbf{SSLPF(10,50,2000)} & 591539.95 & 667071.00 & 749049.35 & 664969.61 & 0.00  & 0.00  & 0.00  & 0.00  & 0.00\% & 0.00\% & 0.00\% & 0.00\% \\
    &       &       &       & (5065.41)  &       &       &       & (0.00)  &       &       &       & (0.00\%) \\
    \hline
    \textbf{SSLPF-indx(10,50,2000)} & 591539.95 & 667071.00 & 749049.35 & 664969.61 & 0.00  & 0.00  & 0.00  & 0.00  & 0.00\% & 0.00\% & 0.00\% & 0.00\% \\
    &       &       &       & (5065.41) &       &       &       & (0.00)  &       &       &       & (0.00\%) \\
    \hline
    \textbf{SSLPF(15,45,15)}       & 22277.30 & 25686.50 & 34244.55 & 26816.72 & 0.00  & 0.00  & 0.00  & 0.00  & 0.00\% & 0.00\% & 0.00\% & 0.00\% \\
    &       &       &       & (480.28) &       &       &       & (0.00)  &       &       &       & (0.00\%) \\
    \hline
    \textbf{SSLPF(15,45,150)}        & 142335.25 & 172721.50 & 223728.10 & 177046.35 & 0.00  & 0.00  & 0.00  & 0.00  & 0.00\% & 0.00\% & 0.00\% & 0.00\% \\
    &       &       &       & (2736.33) &       &       &       & (0.00)  &       &       &       & (0.00\%) \\
    \hline
    \textbf{SSLPF(15,80,15)}    & 180629.90 & 209324.00 & 250452.65 & 210419.66 & 0.00  & 0.00  & 0.00  & 0.00  & 0.00\% & 0.00\% & 0.00\% & 0.00\% \\
    &       &       &       & (2140.48) &       &       &       & (0.00)  &       &       &       & (0.00\%) \\
    \hline   \hline
    \textbf{SMKPF(29)} &       104.35 & 238.50 & 894.25 & 317.21 & 41.00 & 74.50 & 281.3 & 103.47 & 9.78\% & 31.65\% & 71.10\% & 42.08\% \\
    &        &       &       & (24.37) &       &       &       & (16.11)  &       &       &       & (7.65\%) \\
    \hline
    \textbf{SMKPF(30)} &      236.10 & 535.50 & 1947.50 & 783.28 & 45.05 & 78.50 & 381.55 & 147.17 & 4.30\% & 14.03\% & 51.11\% & 20.38\% \\
    &       &       &       & (78.32) &       &       &       & (35.82) &       &       &       & (2.30\%) \\
    \hline
    \multicolumn{13}{l}{\small Standard error of estimate is reported between parentheses.} \\
    \end{tabular}    
    }
    \caption{Relaxed second-stage problems}
    \label{tab:relaxed_second_stage_problems}
\end{table}

\paragraph{Two-phase variants.}
\label{par:feedback_from_approximate_to_exact_solutions}

We now direct our attention to variants of the solution process where ML-L-Shaped outputs a feasible, approximate solution that is used as a warm-start incumbent first-stage solution in Alt-L. We report total computing times in Table~\ref{tab:total_computation_times_with_feedback_from_approximate_solution} along with the average ratio between the two-phase variants of ML-L-Shaped and Alt-L. 
Since this average ratio varies between 83\% and 111\% in the upper panel of the table, we conclude that the application of this two-phase mechanism is not advantageous by itself.

\begin{table}[htbp]
  \centering
  \resizebox{\columnwidth}{!}{
    \begin{tabular}{||l||rrr|r||rrr|r||rrr|r||}
    \hline
    \multirow{3}{*}{\textbf{Problem family}}
    & \multicolumn{4}{c||}{Alt-L} & \multicolumn{4}{c||}{Two-phase and bound} & \multicolumn{4}{c||}{(Two-phase and bound)/Alt-L ratio}    \\
	& \multicolumn{3}{c|}{Quantiles} &  &  \multicolumn{3}{c|}{Quantiles}  &   & \multicolumn{3}{c|}{Quantiles} &    \\
	&  0.05 & 0.5 & 0.95 & Avg & 0.05 & 0.5 & 0.95 & Avg & 0.05 & 0.5 & 0.95 & Avg \\
	\hline
    \hline
    \multicolumn{13}{||c||}{} \\
    \multicolumn{13}{||c||}{\textbf{With warm start only}} \\
    \hline
     \textbf{SSLPF(10,50,2000)} & 131.63 & 150.73 & 186.04 & 156.06 & 113.54 & 127.05 & 143.25 & 127.84 & 72.97\% & 84.36\% & 91.93\% & 82.85\% \\
    &       &       &       & (2.89)  &       &       &       & (1.17)  &       &       &       & (0.67\%) \\
    \hline
    \textbf{SSLPF-indx(10,50,2000)} & 131.63 & 150.73 & 186.04 & 156.06 & 125.57 & 141.74 & 162.68 & 142.03 & 79.88\% & 92.63\% & 102.30\% & 92.02\% \\
    &       &       &       & (2.89)  &       &       &       & (1.43)  &       &       &       & (0.82\%) \\
    \hline
    \textbf{SSLPF(15,45,15)}       & 4.28  & 5.00  & 7.12  & 5.25  & 4.40  & 5.04  & 6.38  & 5.37  & 86.52\% & 100.64\% & 109.28\% & 101.13\% \\
    &       &       &       & (0.11)  &       &       &       & (0.23)  &       &       &       & (1.34\%) \\
    \hline 
    \textbf{SSLPF(15,45,150)}        & 27.96 & 34.16 & 43.93 & 34.96 & 32.14 & 38.46 & 50.46 & 38.78 & 99.40\% & 111.71\% & 120.48\% & 111.22\% \\
    &       &       &       & (0.57)  &       &       &       & (0.62)  &       &       &       & (0.84\%) \\
    \hline
    \textbf{SSLPF(15,80,15)}    & 36.20 & 47.24 & 130.70 & 58.55 & 13.12 & 52.88 & 76.42 & 48.96 & 27.15\% & 98.83\% & 142.27\% & 94.48\% \\
    &       &       &       & (3.15)  &       &       &       & (2.49)  &       &       &       & (6.07\%) \\
    \hline
    \multicolumn{13}{||c||}{} \\
    \multicolumn{13}{||c||}{\textbf{With warm start and probabilistic bound}} \\
    \hline
    \textbf{SSLPF(10,50,2000)} & 131.63 & 150.73 & 186.04 & 156.06 & 99.27 & 112.97 & 130.71 & 113.53 & 60.70\% & 74.97\% & 83.93\% & 73.67\% \\
    &       &       &       & (2.89)  &       &       &       & (1.08)  &       &       &       & (0.72\%) \\
    \hline
    \textbf{SSLPF-indx(10,50,2000)} & 131.63 & 150.73 & 186.04 & 156.06 & 107.54 & 125.34 & 150.39 & 126.95 & 67.78\% & 82.59\% & 95.12\% & 82.29\% \\
    &       &       &       & (2.89)  &       &       &       & (1.43)  &       &       &       & (0.87\%) \\
    \hline
    \textbf{SSLPF(15,45,15)}       & 4.28  & 5.00  & 7.12  & 5.25  & 3.52  & 4.02  & 5.31  & 4.32  & 72.82\% & 81.24\% & 88.53\% & 81.62\% \\
    &       &       &       & (0.11)  &       &       &       & (0.15)  &       &       &       & (0.81\%) \\
    \hline
    \textbf{SSLPF(15,45,150)}        & 27.96 & 34.16 & 43.93 & 34.96 & 24.46 & 29.82 & 41.39 & 30.20 & 78.59\% & 86.72\% & 95.08\% & 86.47\% \\
    &       &       &       & (0.57)  &       &       &       & (0.50)  &       &       &       & (0.53\%) \\
    \hline
    \textbf{SSLPF(15,80,15)}    & 36.20 & 47.24 & 130.70 & 58.55 & 13.57 & 30.56 & 55.99 & 35.29 & 18.79\% & 63.75\% & 114.51\% & 66.73\% \\
    &       &       &       & (3.15)  &       &       &       & (1.57)  &       &       &       & (2.98\%) \\
    \hline   \hline
    \textbf{SMKPF(29)} & 151.82 & 742.89 & 3996.73 & 1233.74 & 167.95 & 859.10 & 4184.03 & 1305.54 & 70.16\% & 114.83\% & 159.15\% & 112.98\% \\    
    &       &       &       & (123.09) &       &       &       & (122.34) &       &       &       & (2.40\%) \\
    \hline
    \textbf{SMKPF(30)} & 212.54 & 1459.19 & 6785.80 & 2132.51 & 251.40	& 1623.59 & 7141.88 &	2350.84	& 79.25\% & 115.25\% & 152.57\% & 115.85\% \\    
    &       &       &       & (254.70) &       &       &       & (244.89) &       &       &       & (1.85\%) \\
    \hline
    \multicolumn{13}{l}{\small Standard error of estimate is reported between parentheses.} \\
    \end{tabular}    
    }
    \caption{Total computing times with warm start from approximate solution (seconds)}
    \label{tab:total_computation_times_with_feedback_from_approximate_solution}
\end{table}

The lower panel of the table reports a similarly calculated average ratio when a probabilistic lower bound on the value of the first-stage objective is introduced in addition in the exact solution process. The joint effect of introducing the warm-start incumbent solution and the probabilistic lower bound on the average time ratio is then moderately yet unambiguously advantageous. (Notice that the application of the probabilistic lower bound alone was shown in a distinct experiment not to be advantageous.) 

\subsection{Families of SMKP instances}
\label{sec:families_of_smkp_instances}

We turn to the analysis of the results related to the SMKP families. Those are reported in the rows below those pertaining to SSLP families in Tables~\ref{tab:computation_times} to~\ref{tab:total_computation_times_with_feedback_from_approximate_solution}. (We do not report results for warm-start only in Table~\ref{tab:total_computation_times_with_feedback_from_approximate_solution} given that the performance with warm start \emph{and} probabilistic bound is poor.) We generate all results related to SMKP by applying the ML-Alt-L variant of ML-L-Shaped. That ML-Alt-L dominates ML-Std-L on SMKP instances was established by independent, preliminary experiments, as expected.

The accuracy of ML-Alt-L reported in Table~\ref{tab:first_stage_values} is excellent. The speed-ups reported in Table~\ref{tab:computation_times} are significant but not of the magnitude witnessed in connection with the SSLP families. The latter is unsurprising in view of the much simpler second-stage problems featured by the SMKP families. Indeed, Tables~\ref{tab:integral_second_stage_problems}~ and~\ref{tab:relaxed_second_stage_problems} show in regard to SMKPF(29) and SMKPF(30) that Alt-L solves the integral and relaxed second-stage problems in only about 50 ms and 5 ms on average, respectively.
In general, large reductions in computing time can take place in contexts where the integral second-stage problems are sufficiently hard. We proceed to illustrate this in the context of the SMKP families.

Provided the two following sufficient conditions hold, the difference in computing times between Alt-L and ML-Alt-L can be made arbitrarily large by increasing the complexity of the second-stage problems in SMKPF(29) and SMKPF(30). \emph{First}, the number of nodes defined by ML-Alt-L and Alt-L must be similar or the computing times of Alt-L and ML-Alt-L must be similar when they handle families SMKPF(29) and SMKPF(30). \emph{Second}, the nearly invariant times required to compute approximate solutions of the integral and relaxed problems of second stage with ML must be at most in the order of 50 ms and 5 ms, respectively. That the first condition holds is illustrated in Table~\ref{tab:number_of_first_stage_nodes}. As regards the second condition, it is shown to hold from calculations made from Tables~\ref{tab:integral_second_stage_problems}~and~\ref{tab:relaxed_second_stage_problems}. Crucially, to make the solution yielded by ML-Alt-L worthwhile, we still need to show that it is sufficiently accurate. Table~\ref{tab:first_stage_values} evidences this.
 
In general, for families of problems characterized by a difficult first stage, ML-Alt-L becomes increasingly advantageous in comparison with its exact version when the difficulty of the second-stage problems increases. The gains in online computing times resulting from ML-Alt-L would become arbitrarily large as the complexity of the second-stage problems would grow unboundedly. Practically, this means that the quantiles and averages of the ratios reported in Table~\ref{tab:computation_times} in regard to families SMKPF(29) and SMKPF(30) would decrease as their numbers of scenarios would grow. In other words, the speed-up would increase as the number of scenarios would grow. Clearly, a practical upper limit on the complexity of the second-stage problems is set by the computation of the corresponding exact solutions when generating the data required for training/validating the ML predictors. One particular manner in which this practical upper limit may be offset is the use of the implicit ML predictor discussed in Section~\ref{sec:families_of_SSLP_instances}.

\section{Comparison with an alternative heuristic}
\label{sec:comparison_with_alternative_heuristic}

Whereas Section~\ref{sec:detailed_analysis} compares the performance of ML-L-Shaped to that of its exact counterpart, it is desirable to also draw a parallel between the performance of ML-L-Shaped and that of a best alternative heuristic algorithm. 
As highlighted in Section~\ref{sec:intro}, the PH algorithm \citep{RockWets91, WatsWood11} is a strong candidate for this purpose.
Tables~\ref{tab:PH_first_stage_values} and~\ref{tab:PH_computation_times} compare the performance of the PH algorithm with that of ML-L-Shaped for each family of problems in classes SSLP and SMKP. 

PH computations are performed on the same workstation as that used to produce the results reported in Section~\ref{sec:detailed_analysis} and are implemented with version~3.8 of the Python language using the Pyomo~5.7.3 optimization library and its PH subpackage PySP \citep{WatsEtAl12}. PySP relies on CPLEX~12.10 to solve the scenario-specific optimization problems defined by the PH algorithm. Parallel computations are handled through MPI by the Mpi4py package. We select the specific settings of the PH algorithm in PySP so as to closely match those demonstrated in \citet{Wood16} to yield a better performance in connection with problems in the SSLP class. We verified that the computing times thus obtained with PH for the original problems SSLP(10, 50, 2000) and SSLP(15, 45, 15) are in the vicinity of those reported in \citet{Wood16}.

From Table~\ref{tab:PH_first_stage_values} we conclude that ML-L-Shaped on the one hand and the PH algorithm on the other hand both approximate the first-stage optimal solution values very well for every family of problems considered in classes SSLP and SMKP. For example, whereas the ML-Std-L achieves a gap of 0.006\% over the instances of family SSLPF(10, 50, 2000), the corresponding figure yielded by the PH algorithm is equal to 0.078\%.

\begin{table}[htbp]
  \centering
  \small
  \resizebox{\columnwidth}{!}{
    \begin{tabular}{||l||rrr|r||rrr|r||rrr|r||}
    \hline
    \multirow{3}{*}{\textbf{Problem family}}
    & \multicolumn{4}{c||}{ML-L-Shaped (From Table~\ref{tab:first_stage_values})} & \multicolumn{4}{c||}{PH algorithm}    \\
	&  \multicolumn{3}{c|}{Quantiles}  &   & \multicolumn{3}{c|}{Quantiles} &    \\
	& 0.05 & 0.5 & 0.95 & Avg & 0.05 & 0.5 & 0.95 & Avg \\
	\hline
    \hline
    \textbf{SSLPF(10,50,2000)} & 0.000 & 0.000 & 0.036 & 0.006 & 0.003 & 0.036 & 0.345 & 0.078 \\   
    &       &       &       & (0.003)  &       &       &       & (0.011) \\
    \hline
     \textbf{SSLPF(15,45,15)}       & 0.000 & 0.000 & 0.608 & 0.064 & 0.001 & 0.001 & 0.024 & 0.005 \\    
     &       &       &       & (0.019)  &       &       &       & (0.001) \\
    \hline
    \textbf{SSLPF(15,45,150)}        & 0.000 & 0.000 & 1.021 & 1.943 & 0.000 & 0.033 & 0.163 & 0.053 \\    
    &       &       &       & (1.150)  &       &       &       & (0.006) \\   
    \hline
     \textbf{SSLPF(15,80,15)}       & 0.000 & 0.000 & 0.538 & 0.075 & 0.000 & 0.021 & 0.205 & 0.119   \\
     &       &       &       & (0.018)  &       &       &       & (0.052) \\
    \hline
    \hline 
    \textbf{SMKPF(29)} & 0.000 & 0.000 & 0.050 & 0.009 & 0.085 & 0.211 & 0.396 & 0.223 \\
    &       &       &       & (0.002) &       &       &       & (0.009) \\ 
    \hline
     \textbf{SMKPF(30)} &  0.000 & 0.000 & 0.027 & 0.005 & 0.098 & 0.216 & 0.384 & 0.224 \\
    &       &       &       & (0.001) &       &       &       & (0.008) \\
    \hline
    \multicolumn{9}{l}{\small Standard error of estimate is reported between parentheses.} \\
    \end{tabular}    
    }
    \caption{Comparison of optimality gaps (percentage) between PH and ML-L-Shaped}
    \label{tab:PH_first_stage_values}
\end{table}

From the figures reported in Table~\ref{tab:PH_computation_times}, we conclude that online computing times are significantly smaller for ML-L-Shaped than for PH with respect to instances in the SSLP class. This is particularly manifest when the number of scenarios is high. For example, in regard to SSLPF(10, 50, 2000), the average ratio between computing times is 28,242\%, representing a 282x speed-up for ML-L-Shaped. Comparing the average computing times for the PH algorithm in Table~\ref{tab:PH_computation_times} with those of the exact algorithm in Table~\ref{tab:computation_times}, we note that the former are strictly greater or nearly equal for all four families in the SSLP class.

In contrast, the computing times achieved by the PH algorithm over the families in the SMKP class are smaller than those achieved by ML-Alt-L. PH achieves average ratios of 25.01\% and 16.99\% representing average speed-ups of respectively 4x and 6x compared to ML-Alt-L.

\begin{table}[htbp]
  \centering
  \resizebox{\columnwidth}{!}{
    \begin{tabular}{||l||rrr|r||rrr|r||rrr|r||}
    \hline
    \multirow{3}{*}{\textbf{Problem family}}
    & \multicolumn{4}{c||}{ML-L-Shaped (from Table~\ref{tab:computation_times})} & \multicolumn{4}{c||}{PH algorithm} & \multicolumn{4}{c||}{PH/ML-L-Shaped ratio}    \\
	& \multicolumn{3}{c|}{Quantiles} &  &  \multicolumn{3}{c|}{Quantiles}  &   & \multicolumn{3}{c|}{Quantiles} &    \\
	&  0.05 & 0.5 & 0.95 & Avg & 0.05 & 0.5 & 0.95 & Avg & 0.05 & 0.5 & 0.95 & Avg \\
    \hline
    \hline
    \textbf{SSLPF(10,50,2000)} & 0.52 & 0.64 & 2.59 & 0.93 & 156.43 & 181.70 & 505.67 & 225.44 & 13120\% & 28174\% & 49685\% & 28242\% \\   
    &       &       &       & (0.08)  &       &       &       & (12.29) &       &       &       & (1024\%) \\
    \hline
    \textbf{SSLPF(15,45,15)}       & 0.37  & 0.43  & 0.60  & 0.45  & 17.17 & 19.56 & 44.05 & 24.80 & 3416\% & 4759\% & 7559\% & 5305\% \\    
    &       &       &       & (0.01)  &       &       &       & (2.26)  &       &       &       & (340\%) \\  
    \hline
    \textbf{SSLPF(15,45,150)}        & 0.39 & 0.53 & 0.79 & 0.55 & 27.72 & 31.31 & 55.98 & 36.69 & 3956\% & 6570\% & 10548\% & 7271\% \\    
    &       &       &       & (0.01)  &       &       &       & (2.08)  &       &       &       & (569\%) \\    
    \hline
    \textbf{SSLPF(15,80,15)}       & 2.32 & 3.22 & 13.90 & 4.12 & 19.54 & 29.20 & 188.65 & 48.88 & 464\% & 948\% & 5376\% & 1338\%  \\  
    &       &       &       & (0.29)  &       &       &       & (4.59)  &       &       &       & (129\%) \\    
    \hline
    \hline 
    \textbf{SMKPF(29)} & 28.02 & 109.12 & 510.91 & 174.55 & 17.58 & 20.27 & 26.80 & 21.17 & 3.81\% & 18.94\% & 68.63\% & 25.01\% \\    
    &       &       &       & (18.05) &       &       &       & (0.33)  &       &       &       & (2.77\%) \\   
    \hline
    \textbf{SMKPF(30)} & 37.91 & 194.34 & 935.81 & 329.82 & 17.84 & 20.18 & 30.12 & 22.49 & 2.34\% & 11.08\% & 57.72\% & 16.99\% \\    
    &       &       &       & (55.97) &       &       &       & (1.00)  &       &       &       & (1.89\%) \\
    \hline
    \multicolumn{13}{l}{\small Standard error of estimate is reported between parentheses.} \\
    \end{tabular}    
    }
    \caption{Comparison of computing times (seconds) between PH and ML-L-Shaped}
    \label{tab:PH_computation_times}
\end{table}

However, the computing times of ML-L-Shaped are arguably invariant with respect to the number of second-stage scenarios whereas the computing times of the PH algorithm are monotonically increasing. Therefore, online computing times of ML-L-Shaped will eventually be exceeded by those of the PH algorithm when the number of scenarios occurring in families of the SMKP class reaches a sufficient magnitude.
Still, it must be ascertained that this threshold does not have a size so large that it is rendered practically irrelevant. To gain useful evidence about this, we assessed the impact on computing times with the PH algorithm of increasing to 2,000 the number of scenarios in SMKPF(29) and SMKPF(30) from their base value of 20. In regard to SMKPF(29) the average computing times rises from 21.17 to 484.57 seconds, and for SMKPF(30) it rises from 22.49 to 372.60 seconds. (Table~\ref{tab:PH_times_with_additional_scenarios} in Appendix~B reports detailed results.)  It thus appears that the critical number of scenarios for which the computing time associated with the PH algorithm exceeds that of ML-L-Shaped would be in the order of a few thousand. This number is sufficiently small to be practically consequential.

\section{Conclusion} \label{sec:Conclusion}

We proposed ML-L-Shaped, a matheuristic leveraging the strong capabilities of generic ML approximators to accelerate the online solution of mixed-integer linear two-stage stochastic programming problems. We aimed to solve problems where the second stage is highly demanding computationally. Our core idea was to substitute the solutions of the latter stages of Benders decomposition with fast, yet accurate predictions arising from supervised ML. We verified the practical usefulness of ML-L-Shaped by applying it to problems seminally addressed in \citet{LapoLouv1993}. We conducted an extensive empirical analysis grounded in families of problems derived from the SSLP and SMKP classes examined in \citet{AnguEtAl2016}. The SSLP instances featured relatively hard second-stage problems and hence correspond to the type of problems for which ML-L-Shaped is designed.

Our results appear to substantiate convincingly the usefulness of ML-L-Shaped: compared to an exact approach, predicting the second-stage solutions with ML can bring about large reductions in online computing time while sacrificing small losses in first-stage solution quality. We reported optimality gaps close to zero and average speed-ups between 5x and 192x. We also compared the online performance with that of the PH algorithm. The results showed that ML-L-Shaped and PH reach solutions of comparable quality. The computing times achieved by ML-L-Shaped over the SSLP instances are considerably smaller than those yielded by PH. In contrast, PH achieved smaller computing times over the SMKP instances. We showed, however, that ML-L-Shaped, whose computing times are invariant with respect to the number of scenarios, would display an advantage when the latter is large.

We remark that the PH algorithm could in principle also be equipped with an ML predictor. This predictor would return heuristic solution values at a high speed for the scenario-specific integer programming problems that must currently be solved repeatedly with an exact solver. Since the PH algorithm must obtain one solution per scenario instead of one expected solution over all scenarios as in Std-L and Alt-L, the computing times of the resulting hybrid PH-ML algorithm could not be invariant with respect to the number of scenarios. However, the relation between number of scenarios and computing time could be made considerably flatter in comparison with the original PH algorithm through the introduction of the high-speed ML-predictions.

An important direction for future research consists in drawing the efficient frontier between (i) the total time required by data generation and ML learning, (ii) the accuracy of the resulting heuristic first-stage solutions and (iii) the time required for their computation. This endeavour would require exploring the outcomes of broad ranges of approximators, hyperparameter settings and data set sizes. In this context, it would also seem opportune to explore \emph{predict-and-optimize} approaches performing joint optimization of the ML predictor and the first-stage solution of problem (P).

Finally, we note that the idea of introducing approximate substitutes for the second-stage solutions of two-stage decompositions is also applicable in principle to the latter stages of multi-stage decompositions. ML predictions for the solution of a given stage would then be made conditional on coupling variables originating from previous stages and on free parameters in the current stage.

\section*{Acknowledgments}

This research was funded by the Canadian National Railway Company (CN) Chair in Optimization of Railway Operations at Universit\'e de Montr\'eal and a Collaborative Research and Development Grant from the Natural Sciences and Engineering Research Council of Canada (CRD-477938-14). Computations were made on the supercomputer B\'eluga, managed by Calcul Qu\'ebec and Compute Canada. The operation of this supercomputer is funded by the Canada Foundation for Innovation (CFI), the Minist\`ere de l'\'Economie, de la Science et de l'Innovation du Qu\'ebec (MESI) and the Fonds de recherche du Qu\'ebec - Nature et technologies (FRQ-NT). The research was also partially funded by the ``IVADO Fundamental Research Project Grants'' under the project entitled ``Machine Learning for (Discrete) Optimization''. 

\bibliographystyle{plainnat_custom}

\bibliography{References}

\pagebreak

\section*{Appendix A: Data generation and ML times}

\begin{table}[htbp]
  \centering
    \begin{tabular}{||l|c|c||}
    \hline
	\textbf{Problem family} & IP/LP & Avg time per example \\
    \hline
    \hline
    \textbf{SSLPF(10,50,2000)} & IP & 3.24   \\
    \hline
    \textbf{SSLPF-indx(10,50,2000)} & IP & 0.0405   \\
    \hline
    \textbf{SSLPF(15,45,15)} & IP & 0.0493   \\
    \hline
    \textbf{SSLPF(15,45,150)} & IP & 0.151   \\
    \hline
    \textbf{SSLPF(15,80,15)} & IP & 0.582  \\
    \hline
    \hline
    \textbf{SMKPF(29)} & IP & 0.129   \\
    \hline
    \textbf{SMKPF(29)} & LP & 0.0188   \\
    \hline
    \textbf{SMKPF(30)} & IP & 0.0743   \\
    \hline
    \textbf{SMKPF(30)} & LP & 0.0269   \\
    \hline
    \multicolumn{3}{l}{\small IP, LP: output is solution of integral or relaxed 2nd stage problem.} \\
    \end{tabular}
    \caption{Data generation times (seconds)}
    \label{tab:average_generation_time_per_example}
\end{table}

\begin{table}[htbp]
  \centering
    \begin{tabular}{||l|c|c||}
    \hline
	\textbf{Problem family} & IP/LP & Training + Validation time \\
    \hline
    \hline
    \textbf{SSLPF(10,50,2000)} & IP & 25.17   \\
    \hline
    \textbf{SSLPF-indx(10,50,2000)} & IP & 4.185   \\
    \hline
    \textbf{SSLPF(15,45,15)} & IP & 24.45   \\
    \hline
    \textbf{SSLPF(15,45,150)} & IP & 33.59   \\
    \hline
    \textbf{SSLPF(15,80,15)} & IP & 32.19  \\
    \hline
    \hline
    \textbf{SMKPF(29)} & IP & 28.78   \\
    \hline
    \textbf{SMKPF(29)} & LP & 39.69   \\
    \hline
    \textbf{SMKPF(30)} & IP & 23.62   \\
    \hline
    \textbf{SMKPF(30)} & LP & 40.95  \\
    \hline
    \multicolumn{3}{l}{\small IP, LP: output is solution of integral or relaxed 2nd stage problem.} \\
    \end{tabular}
    \caption{ML times (hours)}
    \label{tab:training_and_validation_times}
\end{table}

\pagebreak

\section*{Appendix B: PH times with additional scenarios}
\label{annex:b}

\begin{table}[htbp]
  \centering
  \begin{small}
    \begin{tabular}{||l||ccc|c||ccc|c||ccc|c||}
    \hline
    \multirow{2}{*}{\textbf{Problem family}}
	& \multicolumn{3}{c|}{Quantiles} &   \\
	&  0.05 & 0.5 & 0.95 & Avg  \\
    \hline
    \hline
    \textbf{SMKPF(29)} & 17.58 & 20.27 & 26.80 & 21.17  \\
   (20 scen.)       &       &       & & (0.33)  \\   
    \hline
    \textbf{SMKPF(29)} & 187.67 & 300.83 & 678.51 & 484.57  \\
   (2000 scen.)       &       &       & & (98.27)  \\   
    \hline
    \textbf{SMKPF(30)} & 17.84 & 20.18 & 30.12 & 22.49 \\
    (20 scen.)       &       &       & & (1.00)  \\
    \hline
    \textbf{SMKPF(30)} & 177.11 & 232.61 & 586.05 & 372.60 \\
    (2000 scen.)       &       &       & & (80.58)  \\
    \hline
    \multicolumn{5}{l}{\small Standard error of estimate is reported between parentheses.} \\
    \end{tabular}    
    \end{small}
    \caption{PH times with additional scenarios (seconds) }
    \label{tab:PH_times_with_additional_scenarios}
\end{table}

\end{document}